\documentclass[12pt, a4paper]{article}
\usepackage[latin2]{inputenc}
\usepackage{amssymb}
\usepackage{paralist}
\usepackage{amsmath, calc}
\usepackage{amsfonts}
\usepackage{amsthm}
\usepackage{fullpage}
\usepackage{enumerate}
\usepackage{paralist}
\usepackage{comment}

\newtheorem{theorem}{Theorem}[section] %
\newtheorem{lemma}[theorem]{Lemma} %
\newtheorem{proposition}[theorem]{Proposition} %
\newtheorem{claim}{Claim} %

\usepackage{cite}
\usepackage{color}

\title{Multiplicative Sidon basis\\[1em]
\large\textit{Dedicated to the memory of Javier Cilleruelo}}

\begin{document}
\author{S\'andor Z. Kiss \thanks{Department of Algebra and Geometry, Institute of Mathematics, Budapest University of Technology and Economics, M\H{u}egyetem rkp. 3.,
        H-1111, Budapest, Hungary. \newline
        HUN-REN-BME Stochastics Research Group, M\H{u}egyetem rkp. 3., H-1111 Budapest, Hungary.
        Email: kiss.sandor@ttk.bme.hu. ORCID: 0000-0001-9912-8474
        This author was supported by the NKFIH Grants No. K129335, K146387, KKP 144059.}, Csaba
    S\'andor \thanks{Department of Stochastics, Institute of Mathematics, Budapest University of
        Technology and Economics, M\H{u}egyetem rkp. 3., H-1111, Budapest, Hungary. \newline HUN-REN Alfr\'ed R\'enyi Institute of Mathematics, Re\'altanoda utca 13--15.,
        H-1053 Budapest, Hungary.
        MTA--HUN-REN RI Lend\"ulet ``Momentum'' Arithmetic Combinatorics Research Group, Re\'altanoda utca 13--15., H-1053 Budapest, Hungary. Email:
        sandor.csaba@ttk.bme.hu. ORCID:
        0000-0003-3449-8241
        This author was supported by the NKFIH Grants No. K129335, K146387, KKP 144059.}
}
\date{}
\maketitle
\begin{abstract}
\noindent
    For a $h \ge 2$ integer a sequence of positive integers $A$ is called
    asymptotic basis of order $h$ if every sufficiently large positive integer
    $n$ can be written as a sum of $h$ elements of $A$. We say a sequence of
    positive integers $A$ is a multiplicative Sidon sequence if all the products
    of two elements of the sequence $A$ are distinct. In this paper we prove the
    existence of a multiplicative Sidon set which is an asymptotic basis of
    order 2 by using the probabilistic method.
    
\end{abstract}

\section{Introduction}
A sequence of
positive integers $A$ is called an additive Sidon sequence if all the sums of
two elements of the sequence $A$ are distinct.
For a $h \ge 2$ integer a sequence of positive integers $A$ is called
an asymptotic basis of order $h$ if every sufficiently large positive integer
$n$ can be written as a sum of $h$ elements of $A$. We say a sequence of
positive integers $A$ is a multiplicative Sidon sequence if all the products
of two elements of the sequence $A$ are distinct. We say a sequence $A$ is a
multiplicative basis of order $h$ if every positive integer $n$ can be
written as a product of $h$ elements of $A$. We investigate the sets of
natural numbers which satisfy exactly two of the above properties.
Write
\[
A(x) = \sum_{\substack{a\in A\\ a\le x}}1.
\]
It is well-known \cite{DI} that the number of representations of $n$ as a sum of 2 terms from a set $A$ cannot be a constant from a certain point on. On the other hand, if $A$ is an additive Sidon set, then the number of representations of a positive integer $n$ as a sum of 2 terms from $A$ is at most 1.
This implies that there is no additive basis of order 2 which is an additive
Sidon sequence.
Erd\H{o}s, S\'ark\"ozy and S\'os conjectured \cite{ET} the existence of an asymptotic
basis of order $3$ which is an additive Sidon sequence.
This conjecture was proved in \cite{CP}. It is natural to ask whether there exists a set of natural numbers which is both a multiplicative Sidon set and a multiplicative basis of order $h$.   
It is well-known that a multiplicative Sidon set $A$ has at most
$A(x) \ll \frac{x}{\log x}$ elements \cite{EM}. Moreover, for any multiplicative basis of order $h$, Raikov \cite{RA} proved that
\[
\limsup_{x\rightarrow \infty}\frac{A(x)}{\frac{x}{(\log x)^{1-1/h}}} > 0,
\]
which implies that there does not exist a multiplicative Sidon set which is a
multiplicative basis order of $h$.

Furthermore, a multiplicative basis
must contain 1 and the primes so that it has at least
$\gg \frac{x}{\log x}$ elements, which implies that there is no
multiplicative basis which is also an additive Sidon set. The only remaining interesting problem is to study the connection between additive bases and multiplicative Sidon sets.

\begin{theorem}
    There exists a multiplicative Sidon set which is an asymptotic basis of
    order 2.
\end{theorem}

It is easy to see that there exists an infinite set $A$ which is both an additive and a multiplicative Sidon set; $A = \{2^{2^{n}}: n\in \mathbb{Z}^{+}\}$ is suitable. In the finite case Ruzsa \cite{RI} proved that there exists a set
$A\subset \{1,2,\cdots{} ,N\}$ with at least $(1/2 + o(1))\sqrt{N}$
elements as $N\rightarrow \infty$, which is both additive and
multiplicative Sidon set. Finally, it is clear that the set of positive integers is both an additive and multiplicative basis.

Note that the set of primes form a multiplicative Sidon set and the famous
Goldbach conjecture claims that
the multiplicative Sidon set formed by primes is a basis of order 2
for the even positive integers. On the other hand, Goldbach's conjecture
implies that the set of primes is an asymptotic basis of order $3$.

The existence of multiplicative Sidon basis was initially investigated by the 
first author in collaboration with Javier Cilleruelo, but the work was 
discontinued following Javier's illness and subsequent passing.

In this paper we use the following notations. Let $f(x), g(x)$ be real 
functions with nonnegative values. We write $f(x) \ll g(x)$ if there exists a constant $c > 0$ such that $f(x) \le cg(x)$.
Furthermore, we write
$f(x) \asymp g(x)$, if there exist constants $c_{1}, c_{2} > 0$ such that $c_{1}f(x) \le g(x) \le c_{2}f(x)$. Finally, $d(n)$ denotes the number of positive divisors of a positive integer $n$.

\section{Probabilistic tools}

The proof of Theorem 1 is based on the probabilistic method introduced by
Erd\H{o}s and R\'enyi. One can find an excellent description about this method in the Halberstam - Roth book \cite{HL}.
Let $\Omega$ denote the set of strictly increasing sequences of
positive integers. In this paper we denote the probability of an event
$\mathcal{E}$ by $\mathbb{P}(\mathcal{E})$, and the expectation of a random
variable $\zeta$ by $\mathbb{E}(\zeta)$. The following lemma shows the existence of the suitable probability space.

\begin{lemma}
    Let
    \[
    \alpha_{1}, \alpha_{2}, \alpha_{3} \dots{}
    \]
    be real numbers satisfying
    \[
    0 \le \alpha_{n} \le 1 \hspace*{4mm} (n = 1, 2, \dots{}).
    \]
    Then there exists a probability space ($\Omega$, $\mathcal{X}$, $\mathbb{P}$) with the following two properties:
    \begin{itemize}
        \item[(i)] For every natural number $n$, the event $\mathcal{E}^{(n)} = \{A$:
        $A \in \Omega$, $n \in A\}$ is measurable, and
        $\mathbb{P}(\mathcal{E}^{(n)}) = \alpha_{n}$.
        \item[(ii)] The events $\mathcal{E}^{(1)}$, $\mathcal{E}^{(2)}$, ... are independent.
    \end{itemize}
\end{lemma}
See Theorem 13. in \cite{HL}, p. 142.
We denote the characteristic function of the event $\mathcal{E}^{(n)}$ by
\[
t_{(A, n)} = t_{n} =
\left\{
\begin{aligned}
    1 \textnormal{, if } n \in A \\
    0 \textnormal{, if } n \notin A
\end{aligned} \hspace*{3mm}.
\right.
\]

\noindent Furthermore, we denote the number of solutions of
$a_{i_{1}} + a_{i_{2}} + \dots{} + a_{i_{k}} = n$ by $r_{k}(A, n)$,
where
$a_{i_{1}} \in A$, $a_{i_{2}} \in A$, ...,$a_{i_{k}} \in A$, $1 \le a_{i_{1}} 
< a_{i_{2}} < \dots{} < a_{i_{k}} < n$.
Thus
\[
r_{k}(A, n) = \sum_{\overset{(a_{1}, a_{2}, \dots{}, a_{k}) \in
        \mathbb{N}^{k}}{1
        \le a_{1} < \dots{} < a_{k} < n}\atop {a_{1} + a_{2} + \dots{} + a_{k} =
        n}}t_{(A, a_{1})}t_{(A, a_{2})} \dots{} t_{(A, a_{k})}.
\]

It is clear that $r_{k}(A, n)$ is the sum of random
variables. However for $k > 2$ these variables are not independent because the
same $t_{(A, a_{i})}$ may appear in more than one term. To overcome this
problem we need more advanced probabilistic tools.

In the proof of Theorem 1 we apply the method of J. H. Kim and V. H. Vu.
In the next step, we give a short survey of this method;
one can find more details about it in \cite{JK}, \cite{TA}, \cite{VA},
\cite{VU}. Assume that $t_{1}, t_{2}, \dots{} , t_{n}$ are
independent binary (i.e., all $t_{i}$'s are in $\{0,1\}$) random variables.
Consider a polynomial $Y = Y(t_{1}, t_{2}, \dots{} ,
t_{n})$ in $t_{1}, t_{2}, \dots{} ,
t_{n}$ with degree $k$, where the degree of this polynomial is the maximum of
the sum of the exponents of the monomials. A polynomial $Y$ is said to be totally positive if it can be
written in the form $Y = \sum_{i}e_{i}\Gamma_{i}$, where $e_{i} > 0$ and $\Gamma_{i}$ is a product of some $t_{j}$'s for every index $i$.
Given any multi-index
$\underline{\eta} = (\eta_{1}, \dots{} ,\eta_{n}) \in \mathbb{N}^{n}$, the partial derivative $\partial^{\underline{\eta}}(Y)$ of
$Y$ is defined by
\[
\partial^{\underline{\eta}}(Y) = \left(\frac{\partial}{\partial t_{1}}\right)^{\eta_{1}} \cdots{} \left(\frac{\partial}{\partial t_{n}}\right)^{\eta_{n}}Y(t_{1}, t_{2}, \dots{} ,t_{n}),
\]
and denote the order of $\underline{\eta}$ as $|\underline{\eta}| = \eta_{1} + \dots{} + \eta_{n}$. For any order $d \ge 0$, we denote $\mathbb{E}_{d}(Y) = \max_{\eta: |\eta| = d}\mathbb{E}(\partial^{\underline{\eta}}(Y))$. Then $\mathbb{E}_{0}(Y) = \mathbb{E}(Y)$ and $\mathbb{E}_{d}(Y) = 0$ if $d$ exceeds the degree of $Y$.
Define $\mathbb{E}_{\ge d}(Y) =  \max_{d^{'} \ge d}\mathbb{E}_{d^{'}}(Y)$.
We will apply the following inequality due to Kim and Vu.

\begin{lemma}(J. H. Kim - V. H. Vu)
    For every positive integer $k$ and $Y = Y(t_{1}, t_{2}, \dots{} ,t_{n})$ totally
    positive polynomial of degree $k$, where
    the $t_{i}$'s are independent binary random variables, and for any
    $\lambda > 0$ there exists a constant $d_{k} > 0$ depending only on $k$ such
    that
    \[
    \mathbb{P}\left(|Y-\mathbb{E}(Y)| \ge d_{k}\lambda^{k-\frac{1}{2}}\sqrt{\mathbb{E}_{\ge 0}(Y)\mathbb{E}_{\ge 1}(Y)}\right) = O_{k}\left(e^{-\lambda/4+(k-1)\log n}\right).
    \]
\end{lemma}

\noindent See \cite{VU} for the proof. Informally, this theorem asserts that
when all the derivatives of $Y$ are smaller on average than $Y$ itself, and the degree of $Y$ is small, then $Y$ is concentrated
around its mean. Finally we need the Borel - Cantelli lemma:

\begin{lemma}(Borel - Cantelli)
    Let $\{B_{i}\}$ be a sequence of events in a probability space. If
    \[
    \sum_{j=1}^{+\infty}\mathbb{P}(B_{j}) < \infty,
    \]
    \noindent then with probability 1, at most a finite number of the events
    $B_{j}$ can occur.
\end{lemma}
This is Theorem 7. in \cite{HL}, p. 135.

We will use a special case of Chernoff's inequality (Corollary 1.9. in
\cite{AS}), which states that if $t_i$'s are independent Boolean random
variables and $X = t_1 + \dots{} + t_n$, then for any $\delta > 0$ we have
$$
\mathbb{P}\big(|X - \mathbb{E}(X)| \ge \delta \mathbb{E}(X)\big) \le 2e^{-min(\delta^{2}/4,\delta/2)\mathbb{E}(X)}.
$$

\section{Proof of theorem 1.1}

\subsection{Outline of the proof}

We prove the existence of a set $S$ which is both an asymptotic basis of order 2 and a multiplicative Sidon set. 
We define the random set $A$ by
\[
\mathbb{P}(n\in A) = \frac{1}{n^{3/7}}
\]
for every $n\in \mathbb{Z}^+$.

The proof relies on the following four lemmas. First, we prove that the random set $A$ is an additive asymptotic basis of order 2 in the following stronger 
mean.

\begin{lemma}
    There exists a constant $c_{1} > 0$ such that
    \[
    \left|\left\{k: \frac{n}{4} < k < \frac{n}{2}, k, n-k\in A\right\}\right| 
    > c_{1}n^{1/7}
    \]
    with probability 1, for every sufficiently large $n$.
\end{lemma}

Next, we estimate the number of bad elements, namely those that violate the multiplicative Sidon property. 

\begin{lemma}
    With probability 1, there are only finitely many positive integers $n$ for which there exist $a,a^{'},a^{''}\in A$ with $a^{'} < a < a^{''}$ and $a^{2} = a^{'}\cdot a^{''} = n$.
\end{lemma}

\begin{lemma}
    There exists a constant $c_{2} > 0$ such that
    \[
    \left|\left\{k: \frac{n}{4} < k < \frac{n}{2}, \exists l < a < b < n-k, (n-k)l = ab, k,n-k,l,a,b\in A\right\}\right| \le c_{2}n^{17/147}(\log n)^{5}
    \]
    with probability 1, for every sufficiently large $n$.
\end{lemma}

\begin{lemma}
    There exists a constant $c_{3} > 0$ such that
    \[
    \left |\left \{k: \frac{n}{4} < k < \frac{n}{2}, \exists l < a < b < k, kl = ab, k,n-k,l,a,b\in A\right\}\right| \le c_{3}n^{17/147}(\log n)^{5}
    \]
    with probability 1, for every sufficiently large $n$.
\end{lemma}

Finally, we show that removing the bad elements from $A$ preserves the asymptotic basis property.
It follows from Lemma 3.1-3.4 that there exists a set $B$ and a positive integer $N_{0}$ such that, for every $n > N_{0}$, the equation $n = b^{2} = b^{'}b^{''}$ with $b,b^{'},b^{''}\in B$ and $b^{'} < b < b^{''}$ has no solution and
\[
\left|\left\{k: \frac{n}{4} < k < \frac{n}{2}, k, n-k\in A\right\}\right| > c_{1}n^{1/7},
\]
\[
\left|\left\{k: \frac{n}{4} < k < \frac{n}{2}, \exists l < a < b < n-k, (n-k)l = ab, k,n-k,l,a,b\in A\right\}\right| \le c_{2}n^{17/147}(\log n)^{5},
\]
\[
\left|\left\{k: \frac{n}{4} < k < \frac{n}{2}, \exists l < a < b < k, kl = ab, k,n-k,l,a,b\in A\right\}\right| \le c_{3}n^{17/147}(\log n)^{5}.
\]

We get a suitable set $S$ in the following way. We remove all $b\in B$ for which there exist $u,v,w\in B$, $u,v,w<b$ with $bu = vw$. It follows that the set $S$ is a multiplicative Sidon set. Then for every sufficiently large $n$,
\begin{eqnarray*}
    |\{k: k, n-k\in S\}|
    &\ge &\left|\left\{k: \frac{n}{4} < k < \frac{n}{2}, k, n-k\in B\right\}\right|\\
    &-&\left|\left\{k: \frac{n}{4} < k < \frac{n}{2}, \exists l < a < b < n-k, (n-k)l = ab, k,n-k,l,a,b\in B\right\}\right|\\
    &-&\left|\left\{k: \frac{n}{4} < k < \frac{n}{2}, \exists l < a < n-k, (n-k)l = a^{2}, k,n-k,a\in A\right\}\right|\\
    &-&\left|\left\{k: \frac{n}{4} < k < \frac{n}{2}, \exists l < a < b < k, kl = ab, k,n-k,l,a,b\in B\right\}\right|\\
    &-&\left|\left\{k: \frac{n}{4} < k < \frac{n}{2}, \exists l < a < k, kl = a^{2}, k,n-k,l,a\in B\right\}\right|\\
    &\ge & c_{1}n^{1/7} -  c_{2}n^{17/147}(\log n)^{5} - 0 - c_{3}n^{17/147}(\log n)^{5} - 0 > 0,
\end{eqnarray*}
so that $S$ is an asymptotic basis of order $2$.

\subsection{Proof of Lemma 3.1}
Let
\[
R_{n}(A) = \left\{k: \frac{n}{4} < k < \frac{n}{2}, k, n-k\in A\right\}.
\]
\noindent A simple calculation shows that
$$
\mathbb{E}(|R_{n}(A)|) = \sum_{n/4<k<n/2}(k(n-k))^{-3/7}\ge 
\left(\frac{4}{3n}\right)^{3/7}
\sum_{n/4<k<n/2}k^{-3/7} \gg n^{1/7}.
$$
It follows from Chernoff's inequality with $\delta = 1/2$,
that
$$
\mathbb{P}\big(|R_n(A)| \le \frac{1}{2}\mathbb{E}(|R_n(A)|)\big) \le
2e^{-\frac{1}{16}\mathbb{E}(|R_n(A)|)} < 2e^{-cn^{1/7}} < e^{-2\log n} =
\frac{1}{n^{2}}
$$
if $n$ is large enough. By the Borel - Cantelli lemma we have
$$
|R_n(A)| > \frac{1}{2}\mathbb{E}(|R_n(A)|) \gg n^{1/7}
$$
with probability $1$ for every $n$ large enough.

\subsection{Proof of Lemma 3.2}

For a positive integer $n$ let 
$p_{n} = \mathbb{P}(n\in A: \exists a^{'}, a^{''}\in A, a^{'}a^{''} = n^{2}, a^{'} < a^{''})$. Then
\[
p_{n} \le \sum_{\substack{(a^{'}, a^{''}) \\ a^{'}a^{''} = n^{2}}}\frac{1}{n^{3/7}}\cdot \frac{1}{(a^{'})^{3/7}}\cdot \frac{1}{(a^{''})^{3/7}} =
\frac{d(n^{2})}{n^{9/7}}.
\]
It is well known that $d(n^{2}) = n^{o(1)}$, so that $\sum_{n}p_{n}$ is converges, thus by the Borel-Cantelli lemma, with probability 1, there are only finitely
many $n\in A$ for which there exist $a^{'} < a^{''}$ with $a^{'}a^{''} = n^{2}$.

\subsection{Proof of Lemma 3.3}

We split the proof of Lemma 3.3 into three parts by proving the following propositions.

\begin{proposition}
    With probability 1, for all sufficiently large $n$, there are at most two 
    integers $k$ with $\frac{n}{4} < k < \frac{n}{2}$ for which there exist $l < b < a < n-k$ with $n-k,l,a,b\in A$ and $k\in \{l,a,b\}$ such
    that $(n-k)l = ab$.
\end{proposition}

\begin{proposition}
    With probability 1, for all sufficiently large $n$, there are at most four 
    integers $k$ with $\frac{n}{4} < k < \frac{n}{2}$ and $k,n-k\in A$ for which there exists $d\in A$ with $d < n-k$ and $d \neq k$ such that $\textnormal{gcd}(n-k,d) \ge n^{20/21}$.
\end{proposition}

\begin{proposition}
    There exists a constant $c_{4} > 0$ such that for every sufficiently large $n$ with probability 1,
    \[
    \left|\left\{k: \frac{n}{4} < k < \frac{n}{2}, \exists l < a < b < n-k, k,n-k,l,a,b\in A \textnormal{ distinct with } \right.\right.
    \]
    \[
    \left.\left. \textnormal{gcd}(a, n-k)< n^{20/21}, \textnormal{gcd}(b,n-k) < n^{20/21}, (n-k)l = ab\right\}\right| \le c_{4}n^{17/147}(\log n)^{5}.
    \]
\end{proposition}

The above propositions implies Lemma 3.3 in the following way.

\[
\left|\left\{k: \frac{n}{4} < k < \frac{n}{2}, \exists l < b < a < n-k, k,n-k,l,a,b\in A, (n-k)l = ab\right\}\right|
\]
\[
\le \left|\left\{k: \frac{n}{4} < k < \frac{n}{2}, \exists l < b < a < n-k, k,n-k,l,a,b\in A \textnormal{ are distinct with } (n-k)l = ab\right\}\right|
\]
\[
+ \left|\left\{k: \frac{n}{4} < k < \frac{n}{2}, \exists l < b < a < n-k, k,n-k,l,a,b\in A, k\in \{a,b,l\}, (n-k)l = ab\right\}\right|
\]
\[
\le \Big|\Big\{k: \frac{n}{4} < k < \frac{n}{2}, \exists l < b < a < n-k, k,n-l,l,a,b\in A \textnormal{ are distinct with } \textnormal{gcd}(a, n-k)\ge n^{20/21}
\]
\[
\textnormal{ or } \textnormal{gcd}(b,n-k) \ge n^{20/21}, (n-k)l = ab\Big\}\Big|
\]
\[
+ \Big|\Big\{k: \frac{n}{4} < k < \frac{n}{2}, \exists l < b < a < n-k, k,n-k,l,a,b\in A \textnormal{ are distinct with } \textnormal{gcd}(a, n-k)< n^{20/21},
\]
\[
\textnormal{gcd}(b,n-k) < n^{20/21}, (n-k)l = ab\Big\}\Big|
\]
\[
+ \left|\left\{k: \frac{n}{4} < k < \frac{n}{2}, \exists l < b < a < n-k, k,n-k,l,a,b\in A, k\in \{l,a,b\}, (n-k)l = ab\right\}\right|
\]
\[
\le 4 + c_{4}n^{17/147}(\log n)^{5} + 2 \le c_{2}n^{17/147}(\log n)^{5}. 
\]

Now, we prove Proposition 3.5.

\subsubsection{Proof of proposition 3.5}
We define the events $E_{n}$ and $G_{n}$ in the following way.
$$
E_n=\Big(\exists \textnormal{ three } \frac{n}{4}<k<\frac{n}{2},k,n-k\in A  
\textnormal{ such that } \exists a,b,l\in A, l<a<b<n-k,$$ $$ (n-k)l=ab, k\in \{ l,a,b\}\Big) 
$$
$$
\subseteq \Big (\exists \frac{n}{4}<k_1\ne k_2<\frac{n}{2}, k_1,n-k_1,k_2,n-k_2\in A \textnormal{ such that } \exists a_1,b_1,a_2,b_2\in A, k_i<a_i<b_i<n-k_i,$$ $$(n-k_i)k_i=a_ib_i \Big )
$$
$$
\bigcup \Big(\exists \frac{n}{4}<k_1\ne k_2<\frac{n}{2}, k_1,n-k_1,k_2,n-k_2\in A \textnormal{ such that } \exists l_1,l_2,c_1,c_2\in A, l_i<c_i<n-k_{i},l_i<k_i, $$ $$(n-k_i)l_i=k_ic_i\Big)
= F_n\cup G_n
$$
First we show that the event $F_n$ can be covered by the union of six events as follows. In the equations $(n-k_1)k_1=a_1b_1$ and $(n-k_2)k_2=a_2b_2$, the integers
$k_1,k_2,n-k_1,n-k_2$ are distinct, but $a_1$ or $b_1$ may coincide with $k_2, n-k_2, a_2$ or $b_2$ and  $a_2$ or $b_2$ may coincide with $k_1, n-k_1, a_1$ or
$b_1$. Let
$$
U=U_{k_1,k_2,a_1,a_2,b_1,b_2}=$$ $$\left( \{ a_1,b_1\} \cap \{ a_2,b_2\} \right) \cup \left( \{ a_1,b_1\} \cap \{ k_2,n-k_2\} \right) \cup \left( \{ a_2,b_2\} \cap \{
k_1,n-k_1\} \right)
$$
Here $|\{ a_1,b_1\} \cap \{ a_2,b_2\}|\le 1$, otherwise $\{ a_1,b_1\} = \{ a_2,b_2\}$, that is
$$
(n-k_1)k_1=a_1b_1=a_2b_2=(n-k_2)k_2.
$$
It follows that $k_1=k_2$, a contradiction. Similarly, we get that $| \{ a_1,b_1\} \cap \{ k_2,n-k_2\} |\le 1$ and $| \{ a_2,b_2\} \cap \{ k_1,n-k_1\} |\le 1$. It
follows that $|U|\le 3$. We show that $|U|\le 2$. In contrary, let us suppose that $|U|=3$. Then we may assume that there exist $c,d_1,d_2\in A$ such that $\{
a_1,b_1\} =\{ c,d_1\} $ and $\{ a_2,b_2\} =\{ c,d_2\} $ and $d_1=k_2$ or $n-k_2$ and $d_2=k_1$ or $n-k_1$. If $d_1=k_2$ and $d_2=k_1$, then $(n-k_1)k_1=ck_2$ and
$(n-k_2)k_2=ck_1$, where $k_1<k_2$ and $k_2<k_1$, a contradiction. If $d_1=k_2$ and $d_2=n-k_1$, then $(n-k_1)k_1=ck_2$ and $(n-k_2)k_2=c(n-k_1)$, where $k_1<k_2$
and $n-k_1<n-k_2$, a contradiction. If $d_1=n-k_2$ and $d_2=k_1$, then $(n-k_1)k_1=c(n-k_2)$ and $(n-k_2)k_2=ck_1$, where $n-k_2<n-k_1$ and $k_2<k_1$, a
contradiction. If $d_1=n-k_2$ and $d_2=n-k_1$, then $(n-k_1)k_1=c(n-k_2)$ and $(n-k_2)k_2=c(n-k_1)$, where $n-k_2<n-k_1$ and $n-k_1<n-k_2$, a contradiction. It
follows that $|U|\le 2$. Hence
$$
F_n=\bigcup_{j=0}^2 \Big(\exists \frac{n}{4}<k_1\ne k_2<\frac{n}{2} \textnormal{ such that } \exists a_i,b_i, k_i<a_i<b_i<n-k_i,$$ $$ (n-k_i)k_i=a_ib_i, k_i,n-k_i,a_i,b_i\in A, i=1,2, |U|=j
\Big)= \bigcup_{j=0}^2F_n^{(j)}.
$$
Obviously,
$$
F_n^{(0)}=\Big(\exists \frac{n}{4}<k_1\ne k_2<\frac{n}{2} \textnormal{ such that } \exists a_1,b_1,a_2,b_2, k_i<a_i<b_i<n-k_i,
(n-k_i)k_i=a_ib_i,
$$
$$
k_i,n-k_i,a_i,b_i\in A, i=1,2 \hbox{ are distinct }\Big) =E_n^{(1)}
$$
If $|U|=1$, then either $|\{ a_1,b_1\} \cap \{ a_2,b_2\}|= 1$, $|\{ a_1,b_1\} \cap \{ k_2,n-k_2\}|=0$, $|\{ a_2,b_2\} \cap \{ k_1,n-k_1\}|=0$  or we may assume
that  $|\{ a_1,b_1\} \cap \{ a_2,b_2\}|= 0$, $|\{ a_1,b_1\} \cap \{ k_2,n-k_2\}|=0$, $| \{ a_2,b_2\} \cap \{ k_1,n-k_1\} |=1$. Hence we get that
$$
F_n^{(1)}\subseteq
\Big(\exists \frac{n}{4}<k_1\ne k_2<\frac{n}{2} \textnormal{ such that } \exists c,d_1,d_2, \max\{ k_1,k_2\} <c<n-\max\{ k_1,k_2\},
$$
$$
k_1<d_1<n-k_1, k_2<d_2<n-k_2 \textnormal{ such that } (n-k_1)k_1=cd_1, (n-k_2)k_2=cd_2,$$ $$k_1,n-k_1,c,d_1,k_2,n-k_2,d_2\in A \hbox{ are distinct}\Big) 
$$
$$
\bigcup \Big(\exists \frac{n}{4}<k_1\ne k_2<\frac{n}{2} \textnormal{ such that } \exists a_1,b_1,c_2, k_1 <a_1<b_1<n-k_1, k_2<c_2<n-k_2, \textnormal{ such that }
$$
$$
(n-k_1)k_1=a_1b_1, (n-k_2)k_2=(n-k_1)c_2, k_1,n-k_1,a_1,b_1,k_2,n-k_2,c_2\in A \hbox{ are distinct}\Big)
$$
$$
\bigcup \Big(\exists \frac{n}{4}<k_1\ne k_2<\frac{n}{2} \textnormal{ such that } \exists a_1,b_1,c_2, k_1 <a_1<b_1<n-k_1, k_2<c_2<n-k_2 \textnormal{ such that }
$$
$$
(n-k_1)k_1=a_1b_1, (n-k_2)k_2=k_1c_2, k_1,n-k_1,a_1,b_1,k_2,n-k_2,c_2\in A \hbox{ are distinct }\Big)
$$
$$
=\bigcup_{i=2}^4E_n^{(i)}.
$$
If $|U|=2$, then either we may assume $|\{ a_1,b_1\} \cap \{ a_2,b_2\}|= 1$, $|\{ a_1,b_1\} \cap \{ k_2,n-k_2\}|=0$, $|\{ a_2,b_2\} \cap \{ k_1,n-k_1\}|=1$  or
$|\{ a_1,b_1\} \cap \{ a_2,b_2\}|= 0$, $|\{ a_1,b_1\} \cap \{ k_2,n-k_2\}|=1$, $| \{ a_2,b_2\} \cap \{ k_1,n-k_1\} |=1$. We show that the second case is
impossible. Let us suppose that $|\{ a_1,b_1\} \cap \{ a_2,b_2\}|= 0$, $|\{ a_1,b_1\} \cap \{ k_2,n-k_2\}|=1$, $| \{ a_2,b_2\} \cap \{ k_1,n-k_1\} |=1$. Then
$(n-k_1)k_1=k_2c_1$, $(n-k_2)k_2=k_1c_2$ or $(n-k_1)k_1=k_2c_1$, $(n-k_2)k_2=(n-k_1)c_2$ or $(n-k_1)k_1=(n-k_2)c_1$, $(n-k_2)k_2=k_1c_2$ or
$(n-k_1)k_1=(n-k_2)c_1$, $(n-k_2)k_2=(n-k_1)c_2$ for some $k_1<c_1<n-k_1$, $k_2<c_2<n-k_2$, $c_1,c_2\in A$. If $(n-k_1)k_1=k_2c_1$ and $(n-k_2)k_2=k_1c_2$, then
$k_1<k_2$ and $k_2<k_1$, a contradiction. If $(n-k_1)k_1=k_2c_1$ and $(n-k_2)k_2=(n-k_1)c_2$, then $k_1<k_2$ and $n-k_1<n-k_2$, a contradiction. If
$(n-k_1)k_1=(n-k_2)c_1$ and $(n-k_2)k_2=k_1c_2$, then $n-k_2<n-k_1$ and $k_2<k_1$, a contradiction. If $(n-k_1)k_1=(n-k_2)c_1$ and $(n-k_2)k_2=(n-k_1)c_2$, then
$n-k_2<n-k_1$ and $n-k_1<n-k_2$, a contradiction. It follows that if $|U|=2$, then we may assume $|\{ a_1,b_1\} \cap \{ a_2,b_2\}|= 1$, $|\{ a_1,b_1\} \cap \{
k_2,n-k_2\}|=0$, $|\{ a_2,b_2\} \cap \{ k_1,n-k_1\}|=1$. Hence
$$
F_n^{(2)}\subseteq
\Big(\exists \frac{n}{4}<k_1\ne k_2<\frac{n}{2} \textnormal{ such that } \exists c,d_1, \max\{ k_1,k_2\} <c<n-\max\{ k_1,k_2\}, $$ $$ k_1<d_1<n-k_1
\textnormal{ such that } (n-k_1)k_1=cd_1, (n-k_2)k_2=ck_1, k_1,n-k_1,c,d_1,k_2,n-k_2\in A \hbox{ are distinct }\Big) 
$$
$$
\bigcup \Big(\exists \frac{n}{4}<k_1\ne k_2<\frac{n}{2} \textnormal{ such that } \exists c,d_1, \max\{ k_1,k_2\} <c<n-\max\{ k_1,k_2\}, k_1<d_1<n-k_1
$$
$$
\textnormal{ such that } (n-k_1)k_1=cd_1, (n-k_2)k_2=c(n-k_1), k_1,n-k_1,c,d_1,k_2,n-k_2\in A \hbox{ are distinct }\Big)
$$
$$
=\bigcup_{i=5}^6E_n^{(i)}.
$$
The event $G_n$ will be covered by 20 events as follows. Now the equations are $(n-k_1)l_1=k_1c_1$ and $(n-k_2)l_2=k_2c_2$, where $l_i<k_i,c_i<n-k_i$ for $i=1,2$.
We distinguish cases whether $l_1\ne l_2$ or $l_1=l_2$.

\textbf{Case 1.} If $l_1\ne l_2$ and $l_1\ne k_2$ and $l_2\ne k_1$, then the integers $k_1,n-k_1,l_1,k_2,n-k_2,l_2$ are distinct, so for the set
$$
V=V_{k_1,k_2,l_1,l_2,c_1,c_2}= \{ k_1,n-k_1,l_1,c_1,k_2,n-k_2,l_2,c_2\},
$$
we have $|V|=6,7$ or $8$. Clearly,
$$
\Big(\exists \frac{n}{4}<k_1\ne k_2<\frac{n}{2} \textnormal{ such that } \exists c_1,l_1,c_2,l_2, l_1<c_1<n-k_1, l_1<k_1, l_2<c_2< n-k_2, l_2<k_2,$$ $$(n-k_1)l_1=k_1c_1, (n-k_2)l_2=k_2c_2,k_1,n-k_1,c_1,l_1,k_2,n-k_2,c_2,l_2\in A,  |V|=8\Big) 
$$
$$
=\Big(\exists \frac{n}{4}<k_1\ne k_2<\frac{n}{2} \textnormal{ such that } \exists c_1,l_1,c_2,l_2, l_1<c_1<n-k_1,l_{1} < k_{1}, l_2<c_2<n-k_2,l_{2} < k_{2}$$ $$(n-k_1)l_1=k_1c_1 , (n-k_2)l_2=k_2c_2,
k_1,n-k_1,c_1,l_1,k_2,n-k_2,c_2,l_2\in A \hbox{ are distinct }\Big)
$$
$$
= E_n^{(7)}.
$$
If $|V|=7$, then we may assume that $c_2=n-k_1$ or $l_1$ or $k_1$ or $c_1$. Hence
$$
\Big(\exists \frac{n}{4}<k_1\ne k_2<\frac{n}{2} \textnormal{ such that } \exists c_1,l_1,c_2,l_2\in A, l_i<c_i<n-k_i,l_i<k_i, (n-k_i)l_i=k_ic_i, i=1,2,$$ $$ |V|=7\Big)
$$
$$
= \Big(\exists \frac{n}{4}<k_1\ne k_2<\frac{n}{2} \textnormal{ such that } \exists c_1,l_1,l_2, l_1<c_1<n-k_1,l_1<k_1, l_2<k_2,n-k_1<n-k_2,$$
$$
(n-k_1)l_1=k_1c_1 , (n-k_2)l_2=k_2(n-k_1), k_1,n-k_1,c_1,l_1,k_2,n-k_2,l_2\in A \hbox{ are distinct }\Big)
$$
$$
\bigcup \Big(\exists \frac{n}{4}<k_1\ne k_2<\frac{n}{2} \textnormal{ such that } \exists c_1,l_1,l_2, l_1<c_1<n-k_1,l_1<k_1, l_2<k_2,l_{2}<l_1<n-k_2,$$
$$
(n-k_1)l_1=k_1c_1 , (n-k_2)l_2=k_2l_1,
k_1,n-k_1,c_1,l_1,k_2,n-k_2,l_2\in A \hbox{ are distinct }\Big) 
$$
$$
\bigcup \Big(\exists \frac{n}{4}<k_1\ne k_2<\frac{n}{2} \textnormal{ such that } \exists c_1,l_1,l_2, l_1<c_1<n-k_1,l_1<k_1, l_2<k_2,
$$
$$
(n-k_1)l_1=k_1c_1 , (n-k_2)l_2=k_2k_1, k_1,n-k_1,l_1,c_1,k_2,n-k_2,l_2 \hbox{ are distinct }\Big) 
$$
$$
\bigcup \Big(\exists \frac{n}{4}<k_1\ne k_2<\frac{n}{2}\textnormal{ such that } \exists c,l_1,l_2\in A, l_1<c<n-k_1,l_1<k_1, l_2<c<n-k_2,l_2<k_2, 
$$
$$
(n-k_1)l_1=k_1c, (n-k_2)l_2=k_2c, k_1,n-k_1,c,l_1,k_2,n-k_2,l_2\in A \hbox{ are distinct}\Big)
$$
$$
=\bigcup_{i=8}^{11}E_n^{(i)}.
$$
If $|V|=6$, then $c_1=n-k_2$ or $k_2$ or $l_2$ and $c_2=n-k_1$ or $k_1$ or $l_1$. If $c_1=n-k_2$ and $c_2=n-k_1$, then $(n-k_1)l_1=k_1(n-k_2)$ and
$(n-k_2)l_2=k_2(n-k_1)$, that is $n-k_2<n-k_1$ and $n-k_1<n-k_2$, a contradiction. If $c_1=l_2$ and $c_2=l_1$, then $(n-k_1)l_1=k_1l_2$ and $(n-k_2)l_2=k_2l_1$.
Hence $l_1<l_2$ and $l_2<l_1$, a contradiction. It follows that
$$
\Big(\exists \frac{n}{4}<k_1\ne k_2<\frac{n}{2}\textnormal{ such that }\exists c_1,l_1,c_2,l_2\in A, l_i<c_i<n-k_i,l_i<k_i, (n-k_i)l_i=k_ic_i,
$$
$$
i=1,2 , |V|=6\Big) =
$$
$$
\Big(\exists \frac{n}{4}<k_1\ne k_2<\frac{n}{2}\textnormal{ such that } \exists l_1,l_2, l_i<k_i, (n-k_1)l_1=k_1(n-k_2) , (n-k_2)l_2=k_2k_1,
$$
$$
k_1,n-k_1,l_1,k_2,n-k_2,l_2\in A \hbox{ are distinct }\Big)
$$
$$
\bigcup \Big(\exists \frac{n}{4}<k_1\ne k_2<\frac{n}{2}\textnormal{ such that } \exists l_1,l_2, l_i<k_i, (n-k_1)l_1=k_1(n-k_2) , (n-k_2)l_2=k_2l_1,
$$
$$
k_1,n-k_1,l_1,k_2,n-k_2,l_2\in A \hbox{ are distinct }\Big) 
$$
$$
\bigcup \Big(\exists \frac{n}{4}<k_1\ne k_2<\frac{n}{2}\textnormal{ such that } \exists l_1,l_2, l_i<k_i, (n-k_1)l_1=k_1k_2, (n-k_2)l_2=k_2(n-k_1),
$$
$$
k_1,n-k_1,l_1,k_2,n-k_2,l_2\in A \hbox{ are distinct }\Big) 
$$
$$
\bigcup \Big(\exists \frac{n}{4}<k_1\ne k_2<\frac{n}{2}\textnormal{ such that } \exists l_1,l_2, l_i<k_i, (n-k_1)l_1=k_1k_2, (n-k_2)l_2=k_2k_1,
$$
$$
k_1,n-k_1,l_1,k_2,n-k_2,l_2\in A \hbox{ are distinct }\Big) 
$$
$$
\bigcup \Big(\exists \frac{n}{4}<k_1\ne k_2<\frac{n}{2}\textnormal{ such that } \exists l_1,l_2, l_i<k_i, (n-k_1)l_1=k_1k_2, (n-k_2)l_2=k_2l_1,
$$
$$
k_1,n-k_1,l_1,k_2,n-k_2,l_2\in A \hbox{ are distinct }\Big) 
$$
$$
\bigcup \Big(\exists \frac{n}{4}<k_1\ne k_2<\frac{n}{2}\textnormal{ such that } \exists l_1,l_2, l_i<k_i, (n-k_1)l_1=k_1l_2, (n-k_2)l_2=k_2(n-k_1),
$$
$$
k_1,n-k_1,l_1,k_2,n-k_2,l_2\in A \hbox{ are distinct }\Big) 
$$
$$
\bigcup \Big(\exists \frac{n}{4}<k_1\ne k_2<\frac{n}{2}\textnormal{ such that } \exists l_1,l_2, l_i<k_i, (n-k_1)l_1=k_1l_2, (n-k_2)l_2=k_2k_1,
$$
$$
k_1,n-k_1,l_1,k_2,n-k_2,l_2\in A \hbox{ are distinct }\Big)
$$
$$
=\bigcup_{i=12}^{18}E_n^{(i)}.
$$
\textbf{Case 2.} If $l_1\ne l_2$ and $l_1=k_2$ and $l_2\ne k_1$, then the integers $k_1,n-k_1,k_2,n-k_2,l_2$ are distinct, so for the set
$$
W=W_{k_1,k_2,c_1,l_2,c_2}= \{ k_1,n-k_1,c_1,k_2,n-k_2,c_2,l_2\},
$$
we have $|W|=5$ or $6$ or $7$. Clearly,
$$
\Big(\exists \frac{n}{4}<k_1\ne k_2<\frac{n}{2} \textnormal{ such that } \exists c_1,c_2,l_2\in A, k_2<c_1,k_1<n-k_1,  l_2<c_2,k_2<n-k_2,$$
$$
(n-k_1)k_2=k_1c_1, (n-k_2)l_2=k_2c_2, |W|=7\Big)
$$
$$
= \Big(\exists \frac{n}{4}<k_1\ne k_2<\frac{n}{2}\textnormal{ such that } \exists c_1,c_2,l_2, k_2<c_1<n-k_1,k_2<k_1, l_2<c_2<n-k_2,l_2<k_2,$$ $$(n-k_1)k_2=k_1c_1 , (n-k_2)l_2=k_2c_2,
k_1,n-k_1,c_1,k_2,n-k_2,c_2,l_2\in A \hbox{ are distinct }\Big)
$$
$$
=E_n^{(19)}.
$$
If $|W|=6$, then $c_1\in \{ n-k_2,l_2,k_2\} $ or $c_2\in \{ n-k_1,c_1,k_1 \}$. We show that $c_1\in \{ n-k_2,l_2,k_2 \}$ is impossible. If $c_1=n-k_2$, then
$(n-k_1)k_2=k_1(n-k_2)$, that is $k_2n-k_1k_2=k_1n-k_1k_2$. It follows that $k_1=k_2$, a contradiction. If $c_1=l_2$, then $(n-k_1)k_2=k_1l_2$ and
$(n-k_2)l_2=k_2c_2$, then $k_2<l_2$ and $l_2<k_2$, a contradiction. If $c_1=k_2$, then $(n-k_1)k_2=k_1k_2$, that is $n-k_1=k_1$, a contradiction. It follows that
$$
\Big(\exists \frac{n}{4}<k_1\ne k_2<\frac{n}{2}\textnormal{ such that } \exists c_1,l_1,c_2,l_2\in A, l_i<c_i<n-k_i,l_i<k_i, (n-k_i)l_i=k_ic_i , |W|=6\Big)
$$
$$
= \Big(\exists \frac{n}{4}<k_1\ne k_2<\frac{n}{2}\textnormal{ such that } \exists c_1,l_2, k_2<c_1<n-k_1,k_2<k_1, l_2<k_2<n-k_2,
$$
$$
(n-k_1)k_2=k_1c_1 , (n-k_2)l_2=k_2(n-k_1), k_1,n-k_1,c_1,k_2,n-k_2,l_2\in A \hbox{ are distinct }\Big)
$$
$$
\bigcup \Big(\exists \frac{n}{4}<k_1\ne k_2<\frac{n}{2}\textnormal{ such that } \exists c_1,l_2, k_2<c_1<n-k_1,k_2<k_1, l_2<k_2<n-k_2,
$$
$$
(n-k_1)k_2=k_1c_1 , (n-k_2)l_2=k_2k_1, k_1,n-k_1,c_1,k_2,n-k_2,l_2\in A \hbox{ are distinct }\Big)
$$
$$
\bigcup \Big(\exists \frac{n}{4}<k_1\ne k_2<\frac{n}{2}\textnormal{ such that } \exists c,l_2, k_2<c<n-k_1,k_2<k_1, l_2<k_2,
$$
$$
(n-k_1)k_2=k_1c , (n-k_2)l_2=k_2c, k_1,n-k_1,c,k_2,n-k_2,l_2\in A \hbox{ are distinct }\Big)
$$
$$
=\cup_{i=20}^{22}E_n^{(i)}.
$$
If $|W|=5$, then $c_1\in \{ n-k_2,l_2,k_2\} $ and $c_2\in \{ n-k_1,k_2,c_1,k_1 \}$, which is impossible.

\textbf{Case 3.} If $l_1\ne l_2$ and $l_1=k_2$ and $l_2=k_1$, then $(n-k_1)k_2=k_1c_1$ and $(n-k_2)k_1=k_2c_2$, that is $k_2<k_1$ and $k_1<k_2$, a contradiction.

\textbf{Case 4.} If $l_1=l_2=l$, then $l\ne k_1,k_2$, where $(n-k_1)l=k_1c_1$ and $(n-k_2)l=k_2c_2$. Let
$$
Z=Z_{l,k_1,c_1,k_2,c_2}=\{ l,k_1,n-k_1,c_1,k_2,n-k_2,c_2\}.
$$
The integers $k_1,n-k_1,k_2,n-k_2,l$ are distinct. It follows that $|Z|=7$ or $6$ or $5$. Clearly,
$$
\Big(\exists \frac{n}{4}<k_1\ne k_2<\frac{n}{2}\textnormal{ such that } \exists l,c_1,c_2\in A, l<c_i<n-k_i,l<k_i, (n-k_i)l=k_ic_i , |Z|=7\Big)
$$
$$
= \Big(\exists \frac{n}{4}<k_1\ne k_2<\frac{n}{2}\textnormal{ such that } \exists l,c_1,c_2, l<c_i<n-k_i, l<k_i,$$ $$(n-k_1)l=k_1c_1 , (n-k_2)l=k_2c_2, k_1,n-k_1,c_1,l,k_2,n-k_2,c_2\in A \hbox{ are distinct }\Big)
$$
$$
=E_n^{(23)}.
$$
If $|Z|=6$, then we may assume that $c_1=n-k_2$ or $k_2$ or $c_2$. The case $c_1=c_2$ implies that $\frac{c_1}{l}=\frac{n-k_1}{k_1}=\frac{n}{k_1}-1$ and $\frac{c_2}{l}=\frac{n-k_2}{k_2}=\frac{n}{k_2}-1$, that is $\frac{n}{k_1}-1=\frac{n}{k_2}-1$, $k_1=k_2$, a contradiction. It follows that
$$
\Big(\exists \frac{n}{4}<k_1\ne k_2<\frac{n}{2}\textnormal{ such that } \exists l,c_1,c_2\in A, l<c_i<n-k_i, l < k_i, (n-k_i)l=k_ic_i, i=1,2, |Z|=6\Big)
$$
$$
= \Big(\exists \frac{n}{4}<k_1\ne k_2<\frac{n}{2}\textnormal{ such that } 
\exists l,c_2, l<k_1,c_2,k_2<n-k_2, (n-k_1)l=k_1(n-k_2), (n-k_2)l=k_2c_2
$$
$$
k_1,n-k_1,l,k_2,n-k_2,c_2\in A \hbox{ are distinct }\Big) 
$$
$$
\bigcup \Big(\exists \frac{n}{4}<k_1\ne k_2<\frac{n}{2}\textnormal{ such that } \exists l,c_2, l<k_1,c_2,k_2<n-k_2, (n-k_1)l=k_1k_2 , (n-k_2)l=k_2c_2
$$
$$
k_1,n-k_1,l,k_2,n-k_2,c_2\in A \hbox{ are distinct }\Big)
$$
$$
=\bigcup_{i=24}^{25}E_n^{(i)}.
$$
If $|Z|=5$, then $c_1\in \{ k_2,n-k_2\} $ and $c_2\in \{ k_1,n-k_1\} $. If $c_1=k_2$ and $c_2=k_1$, then $(n-k_1)l=k_1k_2$ and $(n-k_2)l=k_2k_1$, that is
$k_1=k_2$, a contradiction. If $c_2=n-k_1$ and $c_1=n-k_2$, then $(n-k_1)l=k_1(n-k_2)$ and $(n-k_2)l=k_2(n-k_1)$, that is $n-k_1>n-k_2$ and $n-k_2>n-k_1$, a
contradiction. It follows that
$$
\Big(\exists \frac{n}{4}<k_1\ne k_2<\frac{n}{2}\textnormal{ such that } \exists l,c_1,c_2\in A, l<c_i<n-k_i,l<k_i, (n-k_i)l=k_ic_i ,i=1,2, |Z|=5\Big)
$$
$$
= \Big(\exists \frac{n}{4}<k_1\ne k_2<\frac{n}{2}\textnormal{ such that } \exists l<\min\{ k_1,k_2\} , (n-k_1)l=k_1(n-k_2) , (n-k_2)l=k_2k_1
$$
$$
k_1,n-k_1,l,k_2,n-k_2\in A \hbox{ are distinct }\Big)
$$
$$
= E_n^{(26)}.
$$
It follows that $E_n\subseteq \cup_{i=1}^{26}E_n^{(i)}$. We prove that $P(E_n)=O(\frac{1}{n^{1.01}})$. It is enough to show that $\mathbb{P}(E_n^{(i)})=O(\frac{1}{n^{1.01}})$ for $1\le i\le 26$.

$E_n^{(1)}$: If $k_1$ is fixed, then there are at most $d((n-k_1)k_1)=n^{o(1)}$ choices for the pair $(a_1,b_1)$. Similarly, if $k_2$ is fixed, then there are at most $d((n-k_2)k_2)=n^{o(1)}$ choices for the pair $(a_2,b_2)$. So we have $n^{2+o(1)}$ choices for the 8-tuple
$(k_1,n-k_1,a_1,b_1,k_2,n-k_2,a_2,b_2)$, where $a_1,b_1,a_2,b_2\asymp n$. It follows that
$$
\mathbb{P}(E_n^{(1)})=\sum_{(k_1,n-k_1,a_1,b_1,k_2,n-k_2,a_2,b_2)}\frac{1}{k_1^{3/7}}\frac{1}{(n-k_1)^{3/7}}\frac{1}{a_1^{3/7}}\frac{1}{b_1^{3/7}}\frac{1}{k_2^{3/7}}\frac{1}{(n-k_2)^{3/7}}\frac{1}{a_2^{3/7}}\frac{1}{b_2^{3/7}}=
$$
$$
\frac{n^{2+o(1)}}{n^{24/7}}=O\left(\frac{1}{n^{1.01}}\right).
$$
$E_n^{(2)}$: 
We need the following lemma
\begin{lemma}\label{lem1}
    \begin{enumerate}
        \item Let $m$ and $a$ be positive integers. The number of solutions of the equation $x^2\equiv a \pmod{m}$ is at most  $\sqrt{\frac{8m}{3}}$.
        \item Let $m$ and $N$ be positive integers, $gcd(m,N^2)=d$. The number of solutions of the equation $x^2\equiv N^2 \pmod{m}$ is at most $2^{\omega
            (m)+1}\sqrt{d}$, where $\omega (m)$ is the number of prime divisors of the integer $m$.
    \end{enumerate}
\end{lemma}

\begin{proof}
    \begin{enumerate}
        \item This is a very simple proof, we sketch it and leave the details to the reader. Let $f_{m}(a)$ be the number of solutions of the equation $x^2\equiv a \pmod{m}$. If $m = p_{1}^{\alpha_{1}}\cdots{} p_{r}^{\alpha_{r}}$, then by the Chinese remainder theorem, $f_{m}(a) = f_{p_{1}^{\alpha_{1}}}(a)\cdots{} f_{p_{r}^{\alpha_{r}}}(a)$.
        It is easy to see that $f_{2^{\alpha}}(a)\le \sqrt{2}\cdot \sqrt{2^{\alpha}}$, where $\alpha > 0$, $f_{3^{\beta}}(a)\le \frac{2}{\sqrt{3}}\cdot \sqrt{3^{\beta}}$, where $\beta > 0$, and for $p > 3$, we have $f_{p^{\gamma}}(a)\le \sqrt{p^{\gamma}}$, where $\gamma > 0$. These implies the required estimation.

        \item Let $m=p_1^{\alpha_1}p_2^{\alpha_2}\dots p_s^{\alpha_s}$. Let us denote $f_m(N^2)$ the number of solutions of the equation $x^2\equiv N^2 \pmod{m}$, where
        $gcd(m,N^2)=d$. By Chinese Remainder Theorem $$f_m(N^2)=f_{p_1^{\alpha_1}}(N^2)f_{p_2^{\alpha_2}}(N^2)\dots f_{p_s^{\alpha_s}}(N^2),$$ where
        $f_{p_i^{\alpha_i}}(N^2)\le 2$ if $gcd(p_i,N^2)=1$ and $p_i>2$, $f_{p_i^{\alpha_i}}(N^2)\le 4$ if $gcd(p_i,N^2)=1$ and $p_i=2$. If $N=q_1^{\beta
            _1}q_2^{\beta_2}\dots q_t^{\beta_t}$ and $gcd(p_i,N^2)>1$, then $p_i=q_j$ for some $1\le j\le t$. If $2\beta_j \ge \alpha_i$, then $x^2\equiv N^2 \equiv 0
        \pmod{p_i^{\alpha _i}}$. It follows that $f_{p_i^{\alpha_i}}(N^2)=p_i^{\lfloor \frac{\alpha_i}{2}\rfloor }\le \sqrt{p_i^{\alpha _i}}$ and in the canonical
        form of $d$ the exponent of the prime number $p_i$ is $\alpha_i$. If $2\beta_j < \alpha_i$, then $x=yp_i^{\beta_j}$, where $0\le y< p_i^{\alpha_i -\beta_j}$, $\textnormal{gcd}(y,p_{i})=1$
        and $y^2\equiv q_1^{2\beta _1}\dots q_{j-1}^{2\beta_{j-1}}\dots q_{j+1}^{2\beta_{j+1}}\dots q_t^{2\beta _t} \pmod{p_i^{\alpha _i-2\beta_j}}$. Hence the number
        of solutions of the equation $x^2\equiv N^2 \pmod{p_i^{\alpha _i}}$ is at most $2p_i^{\alpha_i-\beta_j-(\alpha_i-2\beta_j)}=2p_i^{\beta _j}$ if $p_i>2$. For $p_i=2$, the number
        of solutions of the equation $x^2\equiv N^2 \pmod{p_i^{\alpha _i}}$ is at most $4p_i^{\alpha_i-\beta_j-(\alpha_i-2\beta_j)}=4p_i^{\beta _j} = 2^{\beta_j + 2}$. 
        In this case in the canonical form of $d$ the exponent of the prime number $p_i$ is $2\beta_j$, which completes the proof.
    \end{enumerate}
\end{proof}

The equation $(n-k_i)k_i=cd_i$, $i=1,2$ implies that $(2k_i-n)^2\equiv n^2 \pmod{c}$, where $\frac{n}{4}<c<\frac{3n}{4}$. Let $gcd(c,n^2)=d$. By Lemma \ref{lem1} (2), the number of suitable pairs
$(k_1,k_2)$ is $O(4^{\omega (c)}d)$. Given the integers $c,k_1$ and $k_2$, the integers $d_1$ and $d_2$ are uniquely determined. The number of
$\frac{n}{4}<c<\frac{3n}{4}$ satisfying $gcd(c,n^2)=d$ is $O(\frac{n}{d})$. Hence the number of 5-tuples $(c,k_1,k_2,d_1,d_2)$ satisfying
$gcd(c,n^2)=d$ is $O(\frac{n}{d}d\max_{k\le n}4^{\omega (k)})=n^{1+o(1)}$.
Since $d\mid n^2$, there are $n^{o(1)}$ possible values of $d$. Hence the number of suitable 5-tuples $(c,k_1,k_2,d_1,d_2)$ is $n^{1+o(1)}$, where $k_1,n-k_1,c,d_1,k_2,n-k_2,d_2\asymp n$. It follows that
$$
P(E_n^{(2)})\le
\sum_{(k_1,n-k_1,c,d_1,k_2,n-k_2,d_2)}\frac{1}{k_1^{3/7}}\frac{1}{(n-k_1)^{3/7}}\frac{1}{c^{3/7}}\frac{1}{d_1^{3/7}}\frac{1}{k_2^{3/7}}\frac{1}{(n-k_2)^{3/7}}\frac{1}{d_2^{3/7}}=
$$
$$
\frac{n^{1+o(1)}}{n^3}=O\left(\frac{1}{n^{1.01}}\right).
$$

$E_n^{(3)}$: For given $\frac{n}{4}<k_1<\frac{n}{2}$ there are $d(k_1(n-k_1))=n^{o(1)}$ possible choices for the pair $(a_1,b_1)$. The
equation $(n-k_2)k_2=(n-k_1)c_2$ yields $(2k_2-n)^2\equiv n^2 \pmod{n-k_1}$. By Lemma \ref{lem1} (1) there are $O(\sqrt{n})$ possible choices for $k_2$. It follows that the
number of 7-tuples $(k_1,n-k_1,a_1,b_1,k_2,n-k_2,c_2)$ is $n^{1.5+o(1)}$, where $k_1,n-k_1,a_1,b_1,k_2,n-k_2,c_2\asymp n$. Hence
$$
\mathbb{P}(E_n^{(3)})=\sum_{(k_1,n-k_1,a_1,b_1,k_2,n-k_2,c_2)}\frac{1}{k_1^{3/7}}\frac{1}{(n-k_1)^{3/7}}\frac{1}{a_1^{3/7}}\frac{1}{b_1^{3/7}}\frac{1}{k_2^{3/7}}\frac{1}{(n-k_2)^{3/7}}\frac{1}{c_2^{3/7}}=
$$
$$
\frac{n^{1.5+o(1)}}{n^3}=O\left(\frac{1}{n^{1.01}}\right).
$$

$E_n^{(4)}$: As in the previous proof of $E_n^{(3)}$, we obtain that 
$P(E_n^{(4)})=O\left(\frac{1}{n^{1.01}}\right)$.

$E_n^{(5)}$: If $(n-k_1)k_1=cd_1$, then $(2k_1-n)^2\equiv n^2 \pmod{c}$, where $\frac{n}{4}<c<\frac{3n}{4}$. By Lemma \ref{lem1} (1) the number of integers $\frac{n}{4}<k_1<\frac{n}{2}$ with $(2k_1-n)^2\equiv n^2 \pmod{c}$ is $O(\sqrt{n})$. For given integers $c$ and $k_1$, there is at most one $\frac{n}{4}<k_2<\frac{n}{2}$ with $(n-k_2)k_2=ck_1$. It follows that the number of suitable 6-tuples $(k_1,n-k_1,c,d_1,k_2,n-k_2)$ is $O(n^{1.5})$. Hence
$$
\mathbb{P}(E_n^{(5)})\ll \frac{n^{1.5}}{n^{18/7}}=O\left(\frac{1}{n^{1.01}}\right).
$$

$E_n^{(6)}$: Arguing as in the previous proof of $E_n^{(5)}$, we get that $P(E_n^{(6)})=O\left(\frac{1}{n^{1.01}}\right)$.

\begin{lemma}\label{lem2}
    \begin{enumerate}
        \item Let $n\in \mathbb{Z}^+$, $0<a<n$. The number of solutions of the equation $(n-x)a=xy$, where $\frac{n}{4}<x<\frac{n}{2}$ is $n^{o(1)}$ as $n\to \infty$.
        \item Let $n\in \mathbb{Z}^+$, $0<b<n$. The number of solutions of the equation $(n-x)y=bx$, where $\frac{n}{4}<x<\frac{n}{2}$ is $n^{o(1)}$ as $n\to \infty$.
    \end{enumerate}
\end{lemma}

\textbf{Proof.} (1): If $(n-x)a=xy$, then $na=x(a+y)$. It follows that there are $d(na)=n^{o(1)}$ possibile choices for the integer $x$. If $a$ and $x$ is given, then there is at most one suitable $y$.

(2): If $(n-x)y=bx$, then $(n-x)(y+b)=nb$. It follows that there are $d(nb)=n^{o(1)}$ possibile choices for $x$. If $b$ and $x$ are
given, then there is at most one suitable integer $y$.

$E_n^{(7)}$: By Lemma \ref{lem2} (1), for a given $l_1$ there are $n^{o(1)}$ possibile choices for the pair $(k_1,c_1)$, where $c_1\asymp l_1$. Moreover, for a given $l_2$ there are 
$n^{o(1)}$ possibile choices for the pair $(k_2,c_2)$, where $c_2\asymp l_2$. It follows that
$$
\mathbb{P}(E_n^{(7)})\ll
\sum_{l_1<n}\frac{1}{l_1^{3/7}}\left(\sum_{(k_1,c_1)}\frac{1}{k_1^{3/7}}\frac{1}{(n-k_1)^{3/7}}\frac{1}{c_1^{3/7}}\left(\sum_{l_2<n}\frac{1}{l_2^{3/7}}\left(
\sum_{(k_2,c_2)}\frac{1}{k_2^{3/7}}\frac{1}{(n-k_2)^{3/7}}\frac{1}{c_2^{3/7}}\right)\right)\right)
$$
$$
\ll
\sum_{l_1<n}\frac{1}{l_1^{3/7}}\left(\sum_{(k_1,c_1)}\frac{1}{n^{6/7}}\frac{1}{l_1^{3/7}}\left(\sum_{l_2<n}\frac{1}{l_2^{3/7}}\left(\sum_{(k_2,c_2)}\frac{1}{n^{6/7}}\frac{1}{l_2^{3/7}}\right)\right)\right)
$$
$$
\ll \sum_{l_1<n}\frac{1}{l_1^{3/7}}\left(
\sum_{(k_1,c_1)}\frac{1}{n^{6/7}}\frac{1}{l_1^{3/7}}\left(\sum_{l_2<n}\frac{1}{l_2^{6/7}}\frac{1}{n^{6/7-o(1)}}\right)\right)
$$
$$
\ll \sum_{l_1<n}\frac{1}{l_1^{3/7}}\left( \sum_{(k_1,c_1)}\frac{1}{n^{6/7}}\frac{1}{l_1^{3/7}}\frac{1}{n^{5/7-o(1)}}\right)
\ll \sum_{l_1<n}\frac{1}{l_1^{6/7}}\frac{1}{n^{11/7-o(1)}}
\ll \frac{1}{n^{10/7-o(1)}}\ll \frac{1}{n^{1.01}}.
$$


$E_n^{(8)}$: 
By Lemma \ref{lem2} (1), for given $l_1$ there are $n^{o(1)}$ possible choices for the pair $(k_1,c_1)$. Using again Lemma \ref{lem2} (2), for given $k_1$ there are at most
$n^{o(1)}$ possible choices for the pair $(k_2,l_2)$. It follows that for given $l_1$ there are $n^{o(1)}$ possibile choices for the 4-tuple $(k_1,c_1,k_2,l_2)$, where $c_1\asymp l_1$ and
$l_2\asymp n$. Hence
$$
\mathbb{P}(E_n^{(8)})\ll \sum_{l_1<n}\frac{1}{l_1^{3/7}}\left(
\sum_{(k_1,c_1,k_2,l_2)}\frac{1}{k_1^{3/7}}\frac{1}{(n-k_1)^{3/7}}\frac{1}{c_1^{3/7}}\frac{1}{k_2^{3/7}}\frac{1}{(n-k_2)^{3/7}}\frac{1}{l_2^{3/7}}\right)
$$
$$
\ll \sum_{l_1<n}\frac{1}{l_1^{6/7}}\left( \sum_{(k_1,c_1,k_2,l_2)}\frac{1}{n^{15/7}} \right) \ll
\frac{n^{o(1)}}{n^{15/7}}\sum_{l_1<n}\frac{1}{l_1^{6/7}}=\frac{1}{n^{2-o(1)}}\ll \frac{1}{n^{1.01}}.
$$

$E_n^{(9)}$: By Lemma \ref{lem2}, for given $l_1$ there are $n^{o(1)}$ possible choices for the pair $(k_1,c_1)$, and at most $n^{o(1)}$ possible choices for the pair $(k_2,l_2)$. It follows that for given
$l_1$ there are $n^{o(1)}$ possible choices for $(k_1,c_1,l_2,k_2)$, where $c_1\asymp l_1$ and $l_2\asymp l_1$. Hence
$$
\mathbb{P}(E_n^{(9)})\ll \sum_{l_1<n}\frac{1}{l_1^{3/7}}\left(
\sum_{(k_1,c_1,k_2,l_2)}\frac{1}{k_1^{3/7}}\frac{1}{(n-k_1)^{3/7}}\frac{1}{c_1^{3/7}}\frac{1}{k_2^{3/7}}\frac{1}{(n-k_2)^{3/7}}\frac{1}{l_2^{3/7}}\right)
$$
$$
\ll \sum_{l_1<n}\frac{1}{l_1^{9/7}}\left( \sum_{(k_1,c_1,k_2,l_2)}\frac{1}{n^{12/7}}\right) \ll \frac{n^{o(1)}}{n^{12/7}}\sum_{l_1<n}\frac{1}{l_1^{9/7}}\ll
\frac{1}{n^{12/7-o(1)}}\ll \frac{1}{n^{1.01}}.
$$


$E_n^{(10)}$: By using the same argument to $E_n^{(8)}$ we get that $\mathbb{P}(E_n^{(10)})=O(\frac{1}{n^{1.01}})$.

$E_n^{(11)}$: By Lemma \ref{lem2} (2), for given integer $c$ there are $n^{o(1)}$ possible choices for the 4-tuple $(k_1,l_1,k_2,l_2)$, where $l_1,l_2\asymp c$. Hence
$$
\mathbb{P}(E_n^{(11)})\ll \sum_{c<n}\frac{1}{c^{3/7}}\left(
\sum_{(k_1,k_2,l_2)}\frac{1}{l_1^{3/7}}\frac{1}{k_1^{3/7}}\frac{1}{(n-k_1)^{3/7}}\frac{1}{k_2^{3/7}}\frac{1}{(n-k_2)^{3/7}}\frac{1}{l_2^{3/7}}\right)
$$
$$
\ll \sum_{c<n}\frac{1}{c^{9/7}}\frac{n^{o(1)}}{n^{12/7}}\ll \frac{1}{n^{1.01}}.
$$

$E_n^{(12)}$: By Lemma \ref{lem2} (1), for given $l_1$ there are $n^{o(1)}$ possible choices for the pair $(k_1,k_2)$. For the given pair $(k_1,k_2)$, there is at most one possible value of $l_2$. It follows that for
given $l_1$ there are $n^{o(1)}$ possible choices for the 3-tuple $(k_1,k_2,l_2)$, where $l_1,l_2\asymp n$. Hence
$$
\mathbb{P}(E_n^{(12)})\ll \sum_{l_1<n,
    l_1\asymp n}\frac{1}{l_1^{3/7}}\left(
\sum_{(k_1,k_2,l_2)}\frac{1}{k_1^{3/7}}\frac{1}{(n-k_1)^{3/7}}\frac{1}{k_2^{3/7}}\frac{1}{(n-k_2)^{3/7}}\frac{1}{l_2^{3/7}}\right)
$$
$$
\ll \sum_{l_1<n, l_1\asymp n}\frac{1}{l_1^{3/7}}\frac{n^{o(1)}}{n^{15/7}}\ll \frac{1}{n^{11/7-o(1)}}\ll  \frac{1}{n^{1.01}}.
$$

$E_n^{(13)}$: By Lemma \ref{lem2}, for given $l_1$ there are $n^{o(1)}$ possible choices for the pair $(k_1,k_2)$. For a given pair $(k_1,k_2)$, there is at most one possible value of $l_2$. It follows that for
given $l_1$ there are $n^{o(1)}$ possibile choices for the 3-tuple $(k_1,k_2,l_2)$, where $l_2\asymp l_1$. Hence
$$
\mathbb{P}(E_n^{(13)})\ll \sum_{l_1<n, l_1\asymp n}\frac{1}{l_1^{3/7}}\left(
\sum_{(k_1,k_2,l_2)}\frac{1}{k_1^{3/7}}\frac{1}{(n-k_1)^{3/7}}\frac{1}{k_2^{3/7}}\frac{1}{(n-k_2)^{3/7}}\frac{1}{l_2^{3/7}}\right)
$$
$$
\ll \sum_{l_1<n, l_1\asymp n}\frac{1}{l_1^{6/7}}\frac{n^{o(1)}}{n^{12/7}}\ll \frac{1}{n^{11/7-o(1)}}\ll \frac{1}{n^{1.01}}.
$$

$E_n^{(14)}$, $E_n^{(15)}$ and $E_n^{(16)}$: As in the previous proof of $E_n^{(13)}$ we get that $P(E_n^{(i)})=O(\frac{1}{n^{1.01}})$ for $i=14$, $15$ and $16$.

$E_n^{(17)}$: By Lemma \ref{lem2} (1), for given $l_1$ there are $n^{o(1)}$ possible choices for the pair $(k_1,l_2)$. For the given pair $(k_1,l_2)$, there is at most one possible value for $k_2$. It follows that for
given $l_1$ there are $n^{o(1)}$ possibile choices for the 3-tuple $(k_1,k_2,l_2)$, where $l_2\asymp l_1$. Hence
$$
\mathbb{P}(E_n^{(17)})\ll \sum_{l_1<n}\frac{1}{l_1^{3/7}}\left(
\sum_{(k_1,k_2,l_2)}\frac{1}{k_1^{3/7}}\frac{1}{(n-k_1)^{3/7}}\frac{1}{k_2^{3/7}}\frac{1}{(n-k_2)^{3/7}}\frac{1}{l_2^{3/7}}\right)
$$
$$
=\sum_{l_1<n}\frac{1}{l_1^{6/7}}\frac{n^{o(1)}}{n^{12/7}}\ll \frac{1}{n^{11/7-o(1)}}\ll \frac{1}{n^{1.01}}.
$$

$E_n^{(18)}$: As in the previous proof of $E_n^{(17)}$, we get that $\mathbb{P}(E_n^{(18)})=O\left(\frac{1}{n^{1.01}}\right)$.

$E_n^{(19)}$:  By Lemma \ref{lem2} (2), for a given $l_2$ there are $n^{o(1)}$ possible choices for the pair $(k_2,c_2)$. For a given $k_2$, there is at most $n^{o(1)}$ possible choices for the pair $(k_1,c_1)$.
It follows that for a given $l_2$ there are $n^{o(1)}$ possible choices for the 4-tuple $(k_1,c_1,k_2,c_2)$, where $c_2\asymp l_2$ and $c_1\asymp n$. Hence
$$
\mathbb{P}(E_n^{(19)})\ll \sum_{l_2<n}\frac{1}{l_2^{3/7}}\left(
\sum_{(k_1,c_1,k_2,c_2)}\frac{1}{k_1^{3/7}}\frac{1}{(n-k_1)^{3/7}}\frac{1}{k_2^{3/7}}\frac{1}{(n-k_2)^{3/7}}\frac{1}{l_2^{3/7}}\frac{1}{c_1^{3/7}}\frac{1}{c_2^{3/7}}\right)
$$
$$
\ll \sum_{l_2<n}\frac{1}{l_2^{6/7}}\frac{n^{o(1)}}{n^{15/7}}\ll \frac{1}{n^{2-o(1)}}\ll \frac{1}{n^{1.01}}.
$$

$E_n^{(20)}$: By Lemma \ref{lem2} (2), for given $l_2$ there are $n^{o(1)}$ possible choices for the pair $(k_1,k_2)$. For the given pair $(k_1,k_2)$, there is  at most one choice for the value of $c_1$. It follows
that for given $l_2$ there are $n^{o(1)}$ possibile choices for the 3-tuple $(k_1,c_1,k_2)$, where $c_1,l_2\asymp n$. Hence
$$
\mathbb{P}(E_n^{(20)})\ll \sum_{l_2<n, l_2\asymp n}\frac{1}{l_2^{3/7}}\left(
\sum_{(k_1,c_1,k_2)}\frac{1}{k_1^{3/7}}\frac{1}{(n-k_1)^{3/7}}\frac{1}{k_2^{3/7}}\frac{1}{(n-k_2)^{3/7}}\frac{1}{c_1^{3/7}}\right)
$$
$$
\ll \sum_{l_2<n, l_2\asymp n}\frac{1}{l_2^{3/7}}\frac{n^{o(1)}}{n^{15/7}}\ll \frac{n}{n^{18/7-o(1)}}\ll \frac{1}{n^{1.01}}.
$$

$E_n^{(21)}$: As in the previous proof of $E_n^{(20)}$, we get that $\mathbb{P}(E_n^{(21)})=O\left(\frac{1}{n^{1.01}}\right)$.

$E_n^{(22)}$: By Lemma \ref{lem2} (2), for given $c$ there are $n^{o(1)}$ possible choices for the pair $(k_1,k_2)$. For a given 3-tuple $(c,k_1,k_2)$, there is at most one possible value of $l_2$. It follows
that for given integer $c$ there are $n^{o(1)}$ possible choices for the 3-tuple $(k_1,k_2,l_2)$, where $c,l_2\asymp n$. Hence
$$
\mathbb{P}(E_n^{(22)})\ll \sum_{c<n, c\asymp n}\frac{1}{c^{3/7}}\left(
\sum_{(k_1,k_2,l_2)}\frac{1}{k_1^{3/7}}\frac{1}{(n-k_1)^{3/7}}\frac{1}{k_2^{3/7}}\frac{1}{(n-k_2)^{3/7}}\frac{1}{l_2^{3/7}}\right)
$$
$$
\ll \sum_{c<n, c\asymp n}\frac{n^{o(1)}}{n^{18/7}}\ll \frac{1}{n^{11/7-o(1)}}\ll \frac{1}{n^{1.01}}.
$$

$E_n^{(23)}$: For a given $l$, there are $n^{o(1)}$ possible choices for the pair $(k_1,c_1)$ and there are $n^{o(1)}$ possible choices for the pair $(k_2,c_2)$. It follows that for given $l$ there are
$n^{o(1)}$ possible choices for the 4-tuple $(k_1,c_1,k_2,c_2)$, where $c_1, c_2\asymp l$. Hence
$$
\mathbb{P}(E_n^{(23)})\ll \sum_{l<n}\frac{1}{l^{3/7}}\left(
\sum_{(k_1,c_1,k_2,c_2)}\frac{1}{k_1^{3/7}}\frac{1}{(n-k_1)^{3/7}}\frac{1}{k_2^{3/7}}\frac{1}{(n-k_2)^{3/7}}\frac{1}{c_1^{3/7}}\frac{1}{c_2^{3/7}}\right)
$$
$$
\ll \sum_{l<n}\frac{1}{l^{9/7}}\frac{n^{o(1)}}{n^{12/7}}\ll \frac{1}{n^{12/7-o(1)}}\ll \frac{1}{n^{1.01}}.
$$

$E_n^{(24)}$: By Lemma \ref{lem2} (2), for given $l$ there are $n^{o(1)}$ possible choices for the pair $(k_1,k_2)$. For a given pair $(l,k_2)$, there is at most 
one possible value of $c_2$. It follows that for
given $l$ there are $n^{o(1)}$ possible choices for $(k_1,k_2,c_2)$, where $l,c_2\asymp n$. Hence
$$
\mathbb{P}(E_n^{(24)})\ll \sum_{l<n, l\asymp n}\frac{1}{l^{3/7}}\left(
\sum_{(k_1,k_2,c_2)}\frac{1}{k_1^{3/7}}\frac{1}{(n-k_1)^{3/7}}\frac{1}{k_2^{3/7}}\frac{1}{(n-k_2)^{3/7}}\frac{1}{c_2^{3/7}}\right)
$$
$$
\ll \sum_{l<n, l\asymp n}\frac{n^{o(1)}}{n^{18/7}}\ll \frac{1}{n^{11/7-o(1)}}\ll \frac{1}{n^{1.01}}.
$$

$E_n^{(25)}$: As in the previous proof of $E_n^{(24)}$, we get that $\mathbb{P}(E_n^{(25)})=O\left(\frac{1}{n^{1.01}}\right)$.

$E_n^{(26)}$: By Lemma \ref{lem2} (2), for given $l$ there are $n^{o(1)}$ possible choices for the pair $(k_1,k_2)$, where $l\asymp n$. Hence
$$
\mathbb{P}(E_n^{(26)})\ll \sum_{l<n, l\asymp n}\frac{1}{l^{3/7}}\left( \sum_{(k_1,k_2)}\frac{1}{k_1^{3/7}}\frac{1}{(n-k_1)^{3/7}}\frac{1}{k_2^{3/7}}\frac{1}{(n-k_2)^{3/7}}\right)
$$
$$
\ll \sum_{l<n}\frac{n^{o(1)}}{n^{15/7}}\ll \frac{1}{n^{8/7-o(1)}}\ll \frac{1}{n^{1.01}}.
$$

We get that $\mathbb{P}(E_n)\le \sum_{i=1}^{26}\mathbb{P}(E_n^{(i)})=O\left(\frac{1}{n^{1.01}}\right).$ By the Borel-Cantelli lemma, with probability $1$ only finite number of the events $E_n$
occur, which proves Proposition 3.5.

\subsubsection{Proof of Proposition 3.6}

If $\frac{n}{4} < k_{1} < k_{2} < k_{3} < k_{4} < k_{5} < \frac{n}{2}$ are integers such that for $i = 1, \dots{} ,5$ there exist integers $d_{i} < n-k_i$, $d_{i}
\neq k_{i}$ such that $(d_{i}, n-k_{i}) > n^{20/21}$, then either all the $k_{i}$s, $d_{i}$s and $n-k_i$s are distinct or $d_{i} = d_{j}$ or $d_i = k_j$ or $d_i=n-k_j$ for some $i\neq j$. It follows that
\[
E_{n} = \Big(\exists \frac{n}{4} < k_{1} < k_{2} < k_{3} < k_{4} < k_{5} < \frac{n}{2}, d_i < n-k_i, \textnormal{gcd}(d_i, n-k_i)>n^{20/21},
\]
\[
d_1, \dots{} ,d_5, k_1, \dots{} ,k_5, n-k_1, \dots{} ,n-k_{5}\in A\Big)
\]
\[
\subseteq \Big(\exists \frac{n}{4} < k_{1} < k_{2} < k_{3} < k_{4} < k_{5} < \frac{n}{2}, d_i < n-k_i, \textnormal{gcd}(d_i, n-k_i)>n^{20/21},
\]
\[
d_1, \dots{} ,d_5, k_1, \dots{} ,k_5, n-k_1, \dots{} ,n-k_{5}\in A \textnormal{ are distinct }\Big)
\]
\[
\bigcup \Big(\exists \frac{n}{4} < k_{1} < k_{2} < \frac{n}{2}, d < n-k_2, \textnormal{gcd}(d, n-k_1)>n^{20/21},\textnormal{gcd}(d, n-k_2)>n^{20/21},
\]
\[
k_1,n-k_1,d,k_2,n-k_{2}\in A \textnormal{ are distinct }\Big)
\]
\[
\bigcup \Big(\exists \frac{n}{4} < k_{1} < k_{2} < \frac{n}{2}, d_2 < n-k_2, \textnormal{gcd}(n-k_2, n-k_1)>n^{20/21},\textnormal{gcd}(d_2, n-k_2)>n^{20/21},
\]
\[
k_1, k_2, n-k_1, n-k_{2}, d_2\in A \textnormal{ are distinct }\Big)
\]
\[
\bigcup \Big(\exists \frac{n}{4} < k_{1} < k_{2} < \frac{n}{2}, \textnormal{gcd}(n-k_2, n-k_1)>n^{20/21},\textnormal{gcd}(k_1, n-k_2)>n^{20/21},
\]
\[
k_1, k_2, n-k_1, n-k_{2}\in A \textnormal{ are distinct }\Big)
\]
\[
\bigcup \Big(\exists \frac{n}{4} < k_{1}, k_{2} < \frac{n}{2}, d_2 < n-k_2, (k_2, n-k_1)>n^{20/21}, \textnormal{gcd}(d_2, n-k_2)>n^{20/21},
\]
\[
k_1, k_2, n-k_1, n-k_{2}, d_2\in A \textnormal{ are distinct }\Big)
\]
\[
\bigcup \Big(\exists \frac{n}{4} < k_{1}<k_{2} < \frac{n}{2}, \textnormal{gcd}(k_1, n-k_2)>n^{20/21}, \textnormal{gcd}(k_2, n-k_1)>n^{20/21},
\]
\[
k_1, k_2, n-k_1, n-k_{2}\in A \textnormal{ are distinct }\Big)
\]
\[
= \bigcup_{i=1}^{6}E_{n}^{(i)}.
\]

Now we estimate the probabilities $\mathbb{P}(E_{n}^{(i)})$s. We prove that $\mathbb{P}(E_{n}^{(i)})=O\left(\frac{1}{n^{1.01}}\right)$ for every $1\le i\le 6$. We will use the following observation.

\begin{lemma}\label{lem3}
    Let $1\le m\le n$. The number of $\frac{n}{4}<k<\frac{n}{2}$ with $gcd(n-k,m)>n^{1/21}$ is at most $n^{1/21+o(1)}$.
\end{lemma}

\textbf{Proof.} If $m \le n$ fixed and
$d = \textnormal{gcd}(n-k, m)$, then $d\mid m$ so the number of the possible values of $d$ is $d(m)= n^{o(1)}$. If $d$ is fixed where $d = \textnormal{gcd}(n-k, m) >
n^{20/21}$, then $d\mid n-k$. So the number of possible values of $n-k$ is at most $\frac{n}{d} \le n^{1/21}$. It follows that for a fixed $m \le n$
there are at most
\begin{equation}
    n^{1/21+o(1)}
\end{equation}
integers $\frac{n}{4} < k < \frac{n}{2}$ with $\textnormal{gcd}(n-k,m)\ge n^{20/21}$, which completes the proof of the Lemma \ref{lem3}.

$E_n^{(1)}:$ If $1 \le i \le 5$, then the possible values of $k_i$ is at most $n/2$. If $k_i$ is fixed, then by Lemma \ref{lem3} the number of $d_i$s with $\textnormal{gcd}(n-k_i, d_i) >  n^{20/21}$ is at
most $n^{1/21+o(1)}$. Since $d_i >  n^{20/21}$, we have
\[
\mathbb{P}(E_{n}^{(1)}) \le \sum_{\substack{(k_{1},\dots{} ,k_{5})\\ \frac{n}{4} < k_{i} < \frac{n}{2}}}\left(\prod_{i=1}^{5}\frac{1}{k_{i}^{3/7}}\frac{1}{(n-k_{i})^{3/7}}\right)\left(
\sum_{\substack{(d_{1},\dots{} ,d_{5})\\ \textnormal{gcd}(n-k_i, d_i) >  n^{20/21}}}\prod_{i=1}^{5}\frac{1}{d_{i}^{3/7}}\right)\ll \]
\[
\sum_{\substack{(k_{1},\dots{} ,k_{5})\\ \frac{n}{4} < k_{i} < \frac{n}{2}}}\left(\prod_{i=1}^{5}\frac{1}{k_{i}^{3/7}}\frac{1}{(n-k_{i})^{3/7}}\right)\left( \frac{n^{1/21+o(1)}}{\left( n^{20/21} \right) ^{3/7} } \right) ^5 \ll \sum_{\substack{(k_{1},\dots{} ,k_{5})\\ \frac{n}{4} < k_{i} < \frac{n}{2}}}\left( \frac{1}{n^{30/7}}\frac{n^{5/21+o(1)}}{n^{100/49}} \right) \ll
\]
\[
\frac{n^{5}}{n^{895/147-o(1)}}\  \ll \frac{1}{n^{160/147-o(1)}}\ll\frac{1}{n^{1.01}}.
\]

$E_n^{(2)}:$ The possible values of $k_1$ is at most $n/4$. If $k_1$ is fixed, then by Lemma \ref{lem3} the number of $d$s is at most $n^{1/21+o(1)}$. If $d$ is fixed, then by Lemma \ref{lem3} the
number of $k_2$s is at most $n^{1/21+o(1)}$. Clearly, $d>n^{20/21}$. It follows that
\[
\mathbb{P}(E_{n}^{(2)}) \le \sum_{\frac{n}{4} < k_{1} < \frac{n}{2}}\frac{1}{k_{1}^{3/7}}\frac{1}{(n-k_{1})^{3/7}}\left(\sum_{\substack{d < n-k_i \\ \textnormal{gcd}(n-k_i, d) >
        n^{20/21}}}\frac{1}{d^{3/7}}\left(\sum_{\substack{\frac{n}{4} < k_{2} < \frac{n}{2} \\ \textnormal{gcd}(n-k_2, d) >
        n^{20/21}}}\frac{1}{k_{2}^{3/7}}\frac{1}{(n-k_{2})^{3/7}}\right)\right)
\]
\[
\ll \sum_{\frac{n}{4} < k_1 < \frac{n}{2}}\frac{1}{n^{6/7}}\frac{n^{1/21+o(1)}}{(n^{20/21})^{3/7}}\cdot \frac{n^{1/21+o(1)}}{n^{6/7}}
\ll n\cdot \frac{n^{2/21+o(1)}}{n^{12/7+20/49}}  = \frac{1}{n^{151/147-o(1)}}\ll \frac{1}{n^{1.01}}.
\]

$E_n^{(3)}:$ Let $\frac{n}{4} < k_{2} < \frac{n}{2}$ be fixed. Then by Lemma \ref{lem3} the number of $d_2$s and $k_1$s are at most $n^{1/21+o(1)}$. By $d_2 \ge n^{20/21}$, we have
\[
\mathbb{P}(E_{n}^{(3)}) \le \sum_{\frac{n}{4} < k_{2} < \frac{n}{2}}\frac{1}{k_{2}^{3/7}}\frac{1}{(n-k_{2})^{3/7}}\left(\sum_{\substack{\frac{n}{4} < k_{1} < \frac{n}{2}
        \\ \textnormal{gcd}(n-k_2,n-k_1) >  n^{20/21}}}\frac{1}{k_1^{3/7}}\frac{1}{(n-k_{1})^{3/7}}\left(\sum_{\substack{d_2 < n-k_{2} \\ \textnormal{gcd}(d_2,n-k_2) >
        n^{20/21}}}\frac{1}{d_{2}^{3/7}}\right)\right)
\]
\[
\ll \sum_{\frac{n}{4} < k_2 < \frac{n}{2}}\frac{1}{n^{6/7}}\frac{n^{1/21+o(1)}}{n^{6/7}}\frac{n^{1/21+o(1)}}{(n^{20/21})^{3/7}}
\ll n\cdot \frac{n^{2/21+o(1)}}{n^{12/7+20/49}}  \ll \frac{1}{n^{151/147-o(1)}}\ll \frac{1}{n^{1.01}}.
\]

$E_n^{(4)}:$ Let $d_1 = \textnormal{gcd}(n-k_2,n-k_1) >  n^{20/21}$, $d_2 = \textnormal{gcd}(k_1,n-k_2) >  n^{20/21}$, $d = \textnormal{gcd}(d_1,d_2)$. Then
$d\mid d_1\mid n-k_1$, $d\mid d_2\mid k_1$ so that $d\mid n$. It follows that the number of possible values of $d$ is $d(n) = n^{o(1)}$. Since $d_1\mid n-k_2$,
$d_2\mid n-k_2$ so that $\textnormal{lcm}(d_1,d_2)\mid n-k_2$. Thus $d = \frac{d_1d_2}{\textnormal{lcm}(d_1,d_2)} \ge \frac{n^{20/21}\cdot n^{20/21}}{n-k_2} \ge n^{19/21}$.
If $d\mid n$, $d\ge n^{19/21}$ is fixed, then the number of the pairs $(k_1,k_2)$ with $d\mid n-k_1$, $d\mid n-k_2$ is at most $n^{4/21}$. Thus we have
\[
\sum_{\substack{(k_1,k_2)\\(n-k_2,k_1) >  n^{20/21}\\(n-k_2,k_1) >
        n^{20/21}}}\frac{1}{k_{1}^{3/7}}\frac{1}{(n-k_{1})^{3/7}}\frac{1}{k_{2}^{3/7}}\frac{1}{(n-k_{2})^{3/7}} \ll \frac{1}{n^{12/7}}\cdot n^{4/21} =
\frac{1}{n^{32/21}}\ll \frac{1}{n^{1.01}}.
\]

$E_n^{(5)}:$ Let $\frac{n}{4} < k_{2} < \frac{n}{2}$ be fixed. Then by Lemma \ref{lem3} the number of $\frac{n}{4} < k_{1} < \frac{n}{2}$ with $\textnormal{gcd}(k_2,n-k_1) >  n^{20/21}$, and the number
of $d_2 < n-k_2$ with $\textnormal{gcd}(d_2,n-k_2) >  n^{20/21}$ are at most $n^{1/21+o(1)}$. By $d_2 \ge n^{20/21}$, we have
\[
\mathbb{P}(E_{n}^{(5)}) \le \sum_{\frac{n}{4} < k_{2} < \frac{n}{2}}\frac{1}{k_{2}^{3/7}}\frac{1}{(n-k_{2})^{3/7}}\left(\sum_{\substack{\frac{n}{4} < k_{1} < \frac{n}{2}
        \\ \textnormal{gcd}(k_2,n-k_1) >  n^{20/21}}}\frac{1}{k_1^{3/7}}\frac{1}{(n-k_{1})^{3/7}}\left(\sum_{\substack{d_2 < n-k_{2} \\ \textnormal{gcd}(d_2,n-k_2) >
        n^{20/21}}}\frac{1}{d_{2}^{3/7}}\right)\right)
\]
\[
\ll \sum_{\frac{n}{4} < k_2 < \frac{n}{2}}\frac{1}{n^{6/7}}\frac{n^{1/21+o(1)}}{n^{6/7}}\frac{n^{1/21+o(1)}}{(n^{20/21})^{3/7}}
\ll n\cdot \frac{n^{2/21+o(1)}}{n^{12/7+20/49}}  = \frac{1}{n^{151/147-o(1)}}\ll \frac{1}{n^{1.01}}.
\]

$E_n^{(6)}:$ Let $d_1 = \textnormal{gcd}(n-k_1,k_2)$, $d_2 = \textnormal{gcd}(n-k_2,k_1)$, $d = \textnormal{gcd}(d_1,d_2)$. Then
$d_1\mid n-k_1-k_2$, $d_2\mid n-k_2-k_1$ so that $\textnormal{lcm}(d_1,d_2)\mid n-k_1-k_2$. It follows that $\textnormal{lcm}(d_1,d_2) < n$.
Then $d\mid d_1\mid n-k_1$, $d\mid d_2\mid k_1$ so that $d\mid n$. It follows that the number of possible values of $d$ is $d(n) = n^{o(1)}$ and $d = \frac{d_1d_2}{\textnormal{lcm}(d_1,d_2)} \ge \frac{n^{20/21}\cdot
    n^{20/21}}{n} \ge n^{19/21}$.

If $d\mid n$, $d\ge n^{19/21}$ is fixed, then there are at most $n^{2/21}$ possibilities for the integer $k_1$ with $d\mid k_1$, $\frac{n}{4}<k_1<\frac{n}{2}$. For a fixed $k_1$, by Lemma \ref{lem3} and $gcd(k_2,n-k_1) >  n^{20/21}$, there are at most $n^{1/21+o(1)}$
possibilities for the integer $k_2$. Thus we have
\[
\mathbb{P}(E_{n}^{(6)}) \ll \sum_{\substack{\frac{n}{4} < k_{1} <k_2< \frac{n}{2}\\\textnormal{gcd}(n-k_2,k_1) >  n^{20/21}\\qtextnormal{gcd}(n-k_1,k_2) >
        n^{20/21}}}\frac{1}{k_{1}^{3/7}}\frac{1}{(n-k_{1})^{3/7}}\frac{1}{k_{2}^{3/7}}\frac{1}{(n-k_{2})^{3/7}} \ll \frac{1}{n^{12/7}}\cdot n^{3/21+o(1)} 
\]
\[
= \frac{1}{n^{11/7-o(1)}}\ll \frac{1}{n^{1.01}}.
\]

It follows that
\[
\mathbb{P}(E_{n}) \le \sum_{i=1}^{6}\mathbb{P}(E_{n}^{(i)}) \ll \frac{1}{n^{1.01}}.
\]
By the Borel-Cantelli lemma, with probability 1 only finite number of the events $E_n$ can occur, which proves Proposition 3.6.

\subsubsection{Proof of Proposition 3.7}

Let $t_m=1$, if $m\in A$ and $t_m=0$ if $m\notin A$. Let
\[
X_n = \sum_{\frac{n}{4} < k < \frac{n}{2}}t_kt_{n-k}\sum_{\substack{l<n-k\\l\neq k}}t_l\sum_{\substack{d:d\mid (n-k)l\\ l<d,\frac{(n-k)l}{d}<n-k\\ d^{2}\neq
        (n-k)l,\textnormal{gcd}(d,n-k) <  n^{20/21} \\ \textnormal{gcd}\left(\frac{(n-k)l}{d},n-k \right)< n^{20/21}}}t_{d}t_{\frac{(n-k)l}{d}}.
\]
Since
\[
\Big |\Big \{k: \frac{n}{4} < k < \frac{n}{2}, \exists l<a<b<n-k,  
\textnormal{gcd}(a,n-k) < n^{20/21}, \textnormal{gcd}(b,n-k) < n^{20/21},
\]
\[
(n-k)l = ab, k,n-k,a,b,l\in A \textnormal{ are distinct } \Big \}\Big |\le X_n,
\]
it is enough to show that with probability 1, there exists a constant $c_{10} > 0$ such that
\[
X_n \le c_{10}n^{17/147}(\log n)^{5}
\]
for every sufficiently large $n$. Since $d(m) = m^{o(1)}$,
\[
\mathbb{E}(X_{n}) \ll \sum_{\frac{n}{4} < k < \frac{n}{2}}\frac{1}{k^{3/7}}\frac{1}{(n-k)^{3/7}}\sum_{l<n-k}\frac{1}{l^{3/7}}\sum_{d\mid (n-k)l}\frac{1}{d^{3/7}}
\frac{1}{\left(\frac{(n-k)l}{d}\right)^{3/7}}
\]
\[
\ll \sum_{\frac{n}{4} < k < \frac{n}{2}}\frac{1}{n^{6/7}}\sum_{l<n-k}\frac{1}{l^{6/7}}\frac{d((n-k)l)}{n^{3/7}} \ll \sum_{\frac{n}{4} < k <
    \frac{n}{2}}\frac{1}{n^{6/7}}\frac{(n-k)^{1/7}\cdot n^{o(1)}}{n^{3/7}} \ll \frac{1}{n^{1/7-o(1)}}.
\]
We will show that all expected values of the partial derivatives of $X_n$ are at most $c_{11}n^{17/147}$. It follows that $\mathbb{E}_{\ge 0}(X_n)=O(n^{17/147})$ and 
$\mathbb{E}_{\ge 1}(X_n) = O(n^{17/147})$, 
so that by the Kim-Vu inequality with $\lambda = 24\log n$
\[
\mathbb{P}\left(|X_n - \mathbb{E}(X_n)| \ge c_{12}(24\log n)^{5}\sqrt{\mathbb{E}_{\ge 0}(X_n)\mathbb{E}_{\ge 1}(X_{n})}\right)\ll e^{-\frac{24\log n}{4}+(5-1)\log n}\le c_{13}\cdot \frac{1}{n^{2}}
\]
if $c_{12}$ and $c_{13}$ are large enough. By the Borel-Cantelli lemma, we get that with probability 1,
\[
|X_n - \mathbb{E}(X_n)| \le c_{12}(24\log n)^{5}n^{17/147},
\]
if $n$ is large enough, so we have $X_{n} < c_{14}(\log n)^{5}n^{17/147}$, where
$c_{14}$ is a positive constant.

Now we estimate the expectations of the partial derivatives of $X_n$. Obviously,
\[
\mathbb{E}\left(\frac{\partial}{\partial t_{n-k}}X_n\right) \ll \frac{1}{k^{3/7}}\sum_{l<n-k}\frac{1}{l^{3/7}}\left(\sum_{d\mid (n-k)l}
\frac{1}{d^{3/7}}\frac{1}{\left(\frac{(n-k)l}{d}\right)^{3/7}}\right) \ll \frac{1}{n^{3/7}}\sum_{l<n-k}\frac{1}{l^{3/7}}\frac{d((n-k)l)}{((n-k)l)^{3/7}}
\]
\[
\ll \frac{n^{o(1)}}{n^{6/7}}\sum_{l<n-k}\frac{1}{l^{6/7}} \ll \frac{1}{n^{5/7-o(1)}}\ll n^{17/147},
\]
\[
\mathbb{E}\left(\frac{\partial}{\partial t_{k}}X_n\right) \ll \frac{1}{(n-k)^{3/7}}\sum_{l<n-k}\frac{1}{l^{3/7}}\left(\sum_{d\mid (n-k)l}
\frac{1}{d^{3/7}}\frac{1}{\left(\frac{(n-k)l}{d}\right)^{3/7}}\right) \ll \frac{1}{n^{3/7}}\sum_{l<n-k}\frac{1}{l^{3/7}}\frac{d((n-k)l)}{((n-k)l)^{3/7}}
\]
\[
\ll \frac{n^{o(1)}}{n^{6/7}}\sum_{l<n-k}\frac{1}{l^{6/7}} \ll \frac{1}{n^{5/7-o(1)}}\ll n^{17/147},
\]
\[
\mathbb{E}\left(\frac{\partial}{\partial t_{l}}X_n\right) = \sum_{\frac{n}{4} < k < \frac{n}{2}}\frac{1}{k^{3/7}}\frac{1}{(n-k)^{3/7}}\sum_{d\mid (n-k)l}\frac{1}{d^{3/7}}
\frac{1}{\left(\frac{(n-k)l}{d}\right)^{3/7}} \ll \sum_{\frac{n}{4} < k < \frac{n}{2}}\frac{1}{n^{6/7}}\frac{d((n-k)l)}{((n-k)l)^{3/7}}
\]
\[
\ll \sum_{\frac{n}{4} < k < \frac{n}{2}}\frac{n^{o(1)}}{n^{6/7}n^{3/7}l^{3/7}} \ll \frac{1}{n^{2/7-o(1)}}\ll n^{17/147},
\]
\[
\mathbb{E}\left(\frac{\partial}{\partial t_{d}}X_n\right) = \sum_{\frac{n}{4} < k < \frac{n}{2}}\frac{1}{k^{3/7}}\frac{1}{(n-k)^{3/7}}\sum_{\substack{l<n-k \\ d\mid
        (n-k)l}}\frac{1}{l^{3/7}}
\frac{1}{\left(\frac{(n-k)l}{d}\right)^{3/7}} \ll \sum_{\frac{n}{4} < k < \frac{n}{2}}\frac{1}{n^{3/7}}\sum_{\substack{l<n-k \\ d\mid
        (n-k)l}}\frac{d^{3/7}}{((n-k)l)^{6/7}}
\]
\[
\ll \frac{d^{3/7}}{n^{3/7}}\sum_{\frac{n}{4} < k < \frac{n}{2}}\sum_{\substack{l<n-k \\ d\mid (n-k)l}}\frac{1}{((n-k)l)^{6/7}}
\]
\[
\ll \frac{d^{3/7}}{n^{3/7}} \sum_{m\le \frac{n^{2}}{d}}\frac{d(md)}{(md)^{6/7}}\ll \frac{1}{n^{3/7-o(1)}d^{3/7}}\sum_{m\le \frac{n^{2}}{d}}\frac{1}{m^{6/7}} \ll
\frac{1}{n^{3/7-o(1)}d^{3/7}}\left(\frac{n^{2}}{d}\right)^{1/7} 
\]
\[
\ll \frac{1}{n^{1/7-o(1)}d^{4/7}} \ll n^{17/147}.
\]
Since $\{d: d\mid (n-k)l\} = \{\frac{(n-k)l}{d}: d\mid (n-k)l\}$, thus it follows from the previous expected value that 
\[
\mathbb{E}\left(\frac{\partial}{\partial t_{\frac{(n-k)l}{d}}}X_n\right) \ll n^{17/147}.
\]
Obviously,
\[
\mathbb{E}\left(\frac{\partial^{2}}{\partial t_{n-k} \partial t_{k}}X_n\right) \ll \sum_{l<n-k}\frac{1}{l^{3/7}}\left(\sum_{d\mid (n-k)l}
\frac{1}{d^{3/7}}\frac{1}{\left(\frac{(n-k)l}{d}\right)^{3/7}}\right) \ll \sum_{l<n-k}\frac{1}{l^{3/7}}\frac{d((n-k)l)}{((n-k)l)^{3/7}}
\]
\[
\ll \frac{n^{o(1)}}{n^{3/7}}\sum_{l<n-k}\frac{1}{l^{6/7}} \ll \frac{1}{n^{2/7-o(1)}}\ll n^{17/147},
\]

\[
\mathbb{E}\left(\frac{\partial^{2}}{\partial t_{n-k} \partial t_{l}}X_n\right) \ll \frac{1}{k^{3/7}}\sum_{d\mid (n-k)l}
\frac{1}{d^{3/7}}\frac{1}{\left(\frac{(n-k)l}{d}\right)^{3/7}} \ll \frac{1}{n^{3/7}}\frac{d((n-k)l)}{((n-k)l)^{3/7}}
\]
\[
\ll \frac{1}{n^{6/7-o(1)}}\frac{1}{l^{3/7}}\ll n^{17/147},
\]

\[
\mathbb{E}\left(\frac{\partial^{2}}{\partial t_{n-k} \partial t_{d}}X_n\right) \ll \frac{1}{k^{3/7}}\sum_{l<n-k}
\frac{1}{l^{3/7}}\frac{1}{\left(\frac{(n-k)l}{d}\right)^{3/7}} \ll \frac{d^{3/7}}{n^{6/7}}\sum_{l<n-k}\frac{1}{l^{6/7}}
\]
\[
\ll \frac{d^{3/7}}{n^{5/7}}\ll n^{17/147}.
\]

Since $\{d: d\mid (n-k)l\} = \{\frac{(n-k)l}{d}: d\mid (n-k)l\}$, thus it follows from the previous expected value that

\[
\mathbb{E}\left(\frac{\partial^{2}}{\partial t_{n-k} \partial t_{\frac{(n-k)l}{d}}}X_n\right) \ll n^{17/147}.
\]
Obviously,
\[
\mathbb{E}\left(\frac{\partial^{2}}{\partial t_{k} \partial t_{l}}X_n\right) \ll \frac{1}{(n-k)^{3/7}}\sum_{d\mid (n-k)l}
\frac{1}{d^{3/7}}\frac{1}{\left(\frac{(n-k)l}{d}\right)^{3/7}} \ll \frac{1}{n^{3/7}}\frac{d((n-k)l)}{((n-k)l)^{3/7}}
\]
\[
\ll \frac{n^{o(1)}}{n^{6/7}}\frac{1}{l^{3/7}}\ll n^{17/147},
\]

\[
\mathbb{E}\left(\frac{\partial^{2}}{\partial t_{k} \partial t_{d}}X_n\right) \ll \frac{1}{(n-k)^{3/7}}\sum_{l<n-k}
\frac{1}{l^{3/7}}\frac{1}{\left(\frac{(n-k)l}{d}\right)^{3/7}} \ll \frac{1}{n^{3/7}}\frac{d^{3/7}}{n^{3/7}}\sum_{l<n-k}\frac{1}{l^{6/7}}
\]
\[
\ll \frac{d^{3/7}}{n^{5/7}}\ll n^{17/147}.
\]

Since $\{d: d\mid (n-k)l\} = \{\frac{(n-k)l}{d}: d\mid (n-k)l\}$, thus it follows from the previous expected value that

\[
\mathbb{E}\left(\frac{\partial^{2}}{\partial t_k \partial t_{\frac{(n-k)l}{d}}}X_n\right) \ll n^{17/147}.
\]
Obviously,
\[
\mathbb{E}\left(\frac{\partial^{2}}{\partial t_{l}\partial t_{d}}X_n\right) = \sum_{\substack{\frac{n}{4} < k < \frac{n}{2}\\ d\mid
        (n-k)l}}\frac{1}{k^{3/7}}\frac{1}{(n-k)^{3/7}}
\frac{1}{\left(\frac{(n-k)l}{d}\right)^{3/7}} \ll \frac{1}{n^{9/7}}\sum_{\substack{\frac{n}{4} < k < \frac{n}{2}\\ d\mid (n-k)l}}\frac{d^{3/7}}{l^{3/7}}.
\]
Since $d\mid (n-k)l$ if and only if $\frac{d}{\textnormal{gcd}(d,l)}\mid (n-k)\frac{l}{\textnormal{gcd}(d,l)}$ so that $\frac{d}{\textnormal{gcd}(d,l)}\mid (n-k)$, thus the above sum has $\ll \frac{\textnormal{gcd}(d,l)}{d}n$
terms. Then by $l < d$,
\[
\mathbb{E}\left(\frac{\partial^{2}}{\partial t_{l}\partial t_{d}}X_n\right)\ll \frac{\textnormal{gcd}(d,l)}{d}n\frac{1}{n^{9/7}}\frac{d^{3/7}}{l^{3/7}} \ll
\frac{\textnormal{gcd}(d,l)}{d^{4/7}l^{3/7}n^{2/7}} \ll \frac{l}{d^{4/7}l^{3/7}n^{2/7}}\ll \frac{l^{4/7}}{d^{4/7}n^{2/7}} \ll \frac{1}{n^{2/7}}
\]
\[
\ll n^{17/147}.
\]

Since $\{d: d\mid (n-k)l\} = \{\frac{(n-k)l}{d}: d\mid (n-k)l\}$, thus it follows from the previous expected value that
\[
\mathbb{E}\left(\frac{\partial^{2}}{\partial t_{l} \partial t_{\frac{(n-k)l}{d}}}X_n\right) \ll n^{17/147}.
\]
We have
\[
\mathbb{E}\left(\frac{\partial^{2}}{\partial t_{d} \partial t_{\frac{(n-k)l}{d}}}X_n\right) \ll \sum_{\substack{\frac{n}{4} < k^{'} < \frac{n}{2}\\ n-k^{'}\mid
        d\frac{(n-k)l}{d}}}\frac{1}{(k^{'})^{3/7}}\frac{1}{(n-k^{'})^{3/7}}\sum_{l^{'}\mid d\frac{(n-k^{'})l}{d}}\frac{1}{(l^{'})^{3/7}}
\]
\[
\ll \frac{d(d\frac{(n-k)l}{d})}{n^{6/7}}n^{o(1)} \ll \frac{n^{o(1)}}{n^{6/7}}\ll n^{17/147}.
\]
Clearly,
\[
\mathbb{E}\left(\frac{\partial^{3}}{\partial t_{n-k} \partial t_{k} \partial t_{l}}X_n\right) \ll \sum_{d\mid
    (n-k)l}\frac{1}{d^{3/7}}\frac{1}{\left(\frac{(n-k)l}{d}\right)^{3/7}} \ll \frac{d((n-k)l)}{((n-k)l)^{3/7}}
\]
\[
\ll  \frac{n^{o(1)}}{n^{3/7}l^{3/7}}\ll n^{17/147}.
\]
We have
\[
\mathbb{E}\left(\frac{\partial^{3}}{\partial t_{n-k} \partial t_{k} \partial t_{d}}X_n\right) \ll \sum_{\substack{l<n-k\\ d\mid
        (n-k)l}}\frac{1}{l^{3/7}}\frac{1}{\left(\frac{(n-k)l}{d}\right)^{3/7}}.
\]
Now we have
\[
\frac{n-k}{\textnormal{gcd}(n-k,d)}l = \frac{d}{\textnormal{gcd}(n-k,d)}\cdot \frac{(n-k)l}{d},
\]
so that there exists an $u < \textnormal{gcd}(n-k, d)$ such that $l = u\cdot \frac{d}{\textnormal{gcd}(n-k, d)}$. It follows that
\[
\mathbb{E}\left(\frac{\partial^{3}}{\partial t_{n-k} \partial t_{k} \partial t_{d}}X_n\right) \ll \frac{d^{3/7}}{n^{3/7}}\sum_{u<\textnormal{gcd}(n-k, d)}\frac{1}{\left(u\cdot
    \frac{d}{\textnormal{gcd}(n-k,d)}\right)^{6/7}}
\]
\[
\ll \frac{\textnormal{gcd}(n-k,d)^{6/7}}{n^{3/7}d^{3/7}}\sum_{u<\textnormal{gcd}(n-k, d)}\frac{1}{n^{6/7}} \ll \frac{\textnormal{gcd}(n-k,d)}{n^{3/7}d^{3/7}} \ll \frac{\min\{n^{20/21},d\}}{n^{3/7}d^{3/7}} 
\]
\[
=
\min\left\{\frac{n^{11/21}}{d^{3/7}},\frac{d^{4/7}}{n^{3/7}}\right\} \ll n^{17/147}
\]
because $\frac{n^{11/21}}{d^{3/7}}=\frac{d^{4/7}}{n^{3/7}}$ implies $d = n^{20/21}$, so that
\[
\min\left\{\frac{n^{11/21}}{d^{3/7}},\frac{d^{4/7}}{n^{3/7}}\right\} \le \frac{n^{11/21}}{(n^{20/21})^{3/7}} = n^{17/147}.
\]

Since $\{d: d\mid (n-k)l\} = \{\frac{(n-k)l}{d}: d\mid (n-k)l\}$, thus it follows from the previous expected value that
\[
\mathbb{E}\left(\frac{\partial^{3}}{\partial t_{n-k} \partial t_{k} \partial t_{\frac{(n-k)l}{d}}}X_n\right) \ll n^{17/147}.
\]
Clearly,
\[
\mathbb{E}\left(\frac{\partial^{3}}{\partial t_{n-k} \partial t_{l} \partial t_{d}}X_n\right) \ll \frac{1}{k^{3/7}}\frac{1}{\left(\frac{(n-k)l}{d}\right)^{3/7}} \ll
n^{17/147}.
\]

Since $\{d: d\mid (n-k)l\} = \{\frac{(n-k)l}{d}: d\mid (n-k)l\}$, thus it follows from the previous expected value that
\[
\mathbb{E}\left(\frac{\partial^{3}}{\partial t_{n-k} \partial t_{l} \partial t_{\frac{(n-k)l}{d}}}X_n\right) \ll n^{17/147}.
\]

Clearly,
\[
\mathbb{E}\left(\frac{\partial^{3}}{\partial t_{n-k} \partial t_{d} \partial t_{\frac{(n-k)l}{d}}}X_n\right) \ll \frac{1}{k^{3/7}}\frac{1}{l^{3/7}} \ll
n^{17/147}
\]
and
\[
\mathbb{E}\left(\frac{\partial^{3}}{\partial t_{k} \partial t_{l} \partial t_{d}}X_n\right) \ll \frac{1}{(n-k)^{3/7}}\frac{1}{\left(\frac{(n-k)l}{d}\right)^{3/7}} \ll
n^{17/147}.
\]

Since $\{d: d\mid (n-k)l\} = \{\frac{(n-k)l}{d}: d\mid (n-k)l\}$, thus it follows from the previous expected value that
\[
\mathbb{E}\left(\frac{\partial^{3}}{\partial t_{k} \partial t_{l} \partial t_{\frac{(n-k)l}{d}}}X_n\right) \ll n^{17/147}.
\]

If $l,d,\frac{(n-k)l}{d}$ are given, then the equation $n-k = \frac{(n-k)l}{d}\cdot \frac{d}{l}$ implies that $k$ must be unique, so that
\[
\mathbb{E}\left(\frac{\partial^{3}}{\partial t_{l} \partial t_{d} \partial t_{\frac{(n-k)l}{d}}}X_n\right) = \frac{1}{k^{3/7}}\frac{1}{(n-k)^{3/7}} \ll n^{17/147}.
\]
Furthermore,
\[
\mathbb{E}\left(\frac{\partial^{4}}{\partial t_{n-k} \partial t_{k} \partial t_{l} \partial t_{d}}X_n\right) \ll  \frac{1}{\left(\frac{(n-k)l}{d}\right)^{3/7}} \ll
n^{17/147},
\]

\[
\mathbb{E}\left(\frac{\partial^{4}}{\partial t_{n-k} \partial t_{k} \partial t_{l} \partial t_{\frac{(n-k)l}{d}}}X_n\right) \ll \frac{1}{d^{3/7}} \ll n^{17/147},
\]

\[
\mathbb{E}\left(\frac{\partial^{4}}{\partial t_{n-k} \partial t_{k} \partial t_{d} \partial t_{\frac{(n-k)l}{d}}}X_n\right) \ll \frac{1}{l^{3/7}} \ll n^{17/147},
\]

\[
\mathbb{E}\left(\frac{\partial^{4}}{\partial t_{n-k} \partial t_{l} \partial t_{d} \partial t_{\frac{(n-k)l}{d}}}X_n\right) \ll \frac{1}{k^{3/7}} \ll n^{17/147},
\]

\[
\mathbb{E}\left(\frac{\partial^{4}}{\partial t_{k} \partial t_{l} \partial t_{d} \partial t_{\frac{(n-k)l}{d}}}X_n\right) \ll \frac{1}{(n-k)^{3/7}} \ll n^{17/147}.
\]
Finally,
\[
\mathbb{E}\left(\frac{\partial^{5}}{\partial t_{n-k} \partial t_{k} \partial t_{l} \partial t_{d} \partial t_{\frac{(n-k)l}{d}}}X_n\right) \ll 1 \ll n^{17/147},
\]
which proves Proposition 3.7.

\subsection{Proof of Lemma 3.4}

The proof of Lemma 3.4 is based on the following propositions.

\begin{proposition}
    With probability 1, there are at most four $\frac{n}{4} < k < \frac{n}{2}$ such that there exists a $d < k$ with $k,n-k,d\in A$ and $\textnormal{gcd}(d,k) > n^{20/21}$ if $n$
    is large enough.
\end{proposition}

\begin{proof} Clearly,
    \[
    E_n = \left(\exists \frac{n}{4} < k_1 < \dots{} < k_5 < \frac{n}{2}, d_i < k_i ,1\le i\le 5, d_i, k_i, n-k_i\in A, \textnormal{gcd}(d_i,k_i) > n^{20/21}\right)
    \]
    \[
    \subseteq \Big( \exists \frac{n}{4} < k_1 < \dots{} < k_5 < \frac{n}{2}, d_i < k_i, \textnormal{gcd}(d_i,k_i) >
    n^{20/21},1\le i\le 5,
    \]
    \[
    k_1, n-k_1,d_1,k_2, n-k_2,d_2,k_3, n-k_3,d_3,k_4, n-k_4,d_4,k_5, n-k_5,d_5\in A \textnormal{ are distinct}\Big)
    \]
    \[
    \cup \Big(\exists \frac{n}{4} < k_1 < k_2 < \frac{n}{2}, d < k_1,\textnormal{gcd}(d,k_1) > n^{20/21}, \textnormal{gcd}(d,k_2) > n^{20/21}, \]
    \[
    k_1,n-k_1,d, k_2,  n-k_2\in A \textnormal{ are distinct }\Big)
    \]
    \[
    \cup \Big(\exists \frac{n}{4} < k_1 < k_2 < \frac{n}{2}, d_1 < k_1, \textnormal{gcd}(k_1, d_1) > n^{20/21}, \textnormal{gcd}(k_2,k_1) >
    n^{20/21}, \]
    \[
    d_1, k_1, k_2 n-k_1, n-k_2\in A \textnormal{ are distinct }, \Big)
    \]
    \[
    E_{1}^{(n)} \cup E_{2}^{(n)} \cup E_{3}^{(n)}.
    \]
    
    Similary as the estimation of $\mathbb{P}(E_{1}^{(n)})$, $\mathbb{P}(E_{2}^{(n)})$, $\mathbb{P}(E_{3}^{(n)})$ in Proposition 3.6, we get that $$\mathbb{P}(E_{1}^{(n)}) \ll  \frac{1}{n^{160/147-o(1)}},$$
    $$\mathbb{P}(E_{2}^{(n)}) \ll  \frac{1}{n^{151/147-o(1)}}$$, $$\mathbb{P}(E_{3}^{(n)}) \ll  \frac{1}{n^{151/147-o(1)}}.$$ It follows that
    $$\mathbb{P}(E_{n}) \le \mathbb{P}(E_{1}^{(n)}) + \mathbb{P}(E_{1}^{(n)}) + \mathbb{P}(E_{1}^{(n)}) \ll \frac{1}{n^{151/147-o(1)}}.$$
    By the Borel-Cantelli lemma, with probability 1 only finite number of the events $F_n$ occur, which proves Proposition 3.8.
\end{proof}

\begin{proposition}
    There exists a constant $c_{15} > 0$ such that with probability 1,
    \[
    \Big|\Big\{ \frac{n}{4} < k < \frac{n}{2}, \exists l<a<b<k,  kl = ab, k,n-k,a,b,l\in A,
    \textnormal{gcd}(a,k) < n^{20/21},
    \]
    \[
    gcd(b,k) < n^{20/21} \Big\}\Big| < c_{15}n^{17/147}(\log n)^{5}
    \]
    if $n$ is large enough.
\end{proposition}

The proof is very similar to Proposition 3.7. We apply the same idea for
\[
Y_{n} = \sum_{\frac{n}{4} < k < \frac{n}{2}}t_{n-k}t_{k}\sum_{l<k}t_l
\sum_{\substack{l < d < k\\ d\mid kl\\ d^{2}\neq kl}}t_{d}t_{\frac{kl}{d}}.
\]
Similarly to Proposition 3.7, we can prove that all the partial derivatives are $\ll n^{17/147}$. We omit the details.

The proof of Lemma 3.4 comes from these propositions in the following way. Proposition 3.8 and 3.9 imply that with probability 1,
\[
\left|\left\{k:\frac{n}{4} < k < \frac{n}{2}, \exists l<a<b<k,  kl = ab,   k,n-k,a,b,l\in A\right\}\right|
\]
\[
\le \Big |\Big \{k: \frac{n}{4} < k < \frac{n}{2}, \exists l<a<b<k,  kl = ab,
\textnormal{gcd}(a,k) > n^{20/21},\textnormal{ or } \textnormal{gcd}(b,k) > n^{20/21},
\]
\[
k,n-k,a,b,l\in A\Big \}\Big |
\]
\[
+ \Big |\Big \{k: \frac{n}{4} < k < \frac{n}{2}, \exists l<a<b<k,  kl = ab, \textnormal{gcd}(a,k) < n^{20/21},\textnormal{ and } \textnormal{gcd}(b,k) < n^{20/21}, 
\]
\[
k,n-k,a,b,l\in A\Big\}\Big|
\]
\[
\le \Big|\Big\{k: \frac{n}{4} < k < \frac{n}{2}, \exists d<k,  \textnormal{gcd}(d,k) > n^{20/21}, d,k,n-k\in A
\Big\}\Big|
\]
\[
+ \Big|\Big\{k: \frac{n}{4} < k < \frac{n}{2}, \exists l<a<b<k,  kl = ab, \textnormal{gcd}(a,k) < n^{20/21},\textnormal{ and } \textnormal{gcd}(b,k) < n^{20/21}, 
\]
\[
k,n-k,a,b,l\in A
\Big\}\Big|
\]
\[
\le 4 + c_{15}n^{17/147}(\log n)^{5} \le c_{16}n^{14/147}(\log n)^{5},
\]
which completes the proof.

\end{document}